\documentclass[11pt]{amsart}
\usepackage{verbatim}
\usepackage{latexsym}
\usepackage{amsmath,amsthm,amssymb, amscd,color}
\usepackage{mathrsfs}
\usepackage[all]{xy}
\usepackage[all]{xy}
\usepackage{tikz}
\usepackage{tikz-cd}
\usepackage{mathabx}
\usetikzlibrary{decorations.pathmorphing}
\usepackage{anysize}
\marginsize{3cm}{3cm}{2.5cm}{2.5cm}

\usepackage{lmodern,babel,adjustbox,booktabs,multirow}
\usepackage[T1]{fontenc}

\input xy \xyoption{all}
\usepackage{blkarray} 
\usepackage{framed}
\usepackage[utf8]{inputenc}

\theoremstyle{plain}
\newtheorem{theorem}{Theorem}[section]

\numberwithin{equation}{section}

\theoremstyle{definition}
\newtheorem{definition}[theorem]{Definition}
\newtheorem{example}[theorem]{Example}

\newtheorem{proposition}[theorem]{Proposition}
\newtheorem{remark}[theorem]{Remark}

\theoremstyle{remark}

\newcommand{\CC}{\mathbb{C}}

\newcommand{\RR}{\mathbb{R}}
\newcommand{\ZZ}{\mathbb{Z}}

\begin{document}

\title[On cohomology of quasitoric manifolds]{On cohomology of quasitoric manifolds over a vertex cut of a finite product of simplices}


\author[S. Sarkar]{Soumen Sarkar}
\address{Department of Mathematics, Indian Institute of Technology Madras, India}
\email{soumen@iitm.ac.in}

\author[S. Sau]{Subhankar Sau}
\address{Department of Mathematics, The Institute of Mathematical Sciences, Chennai, India}
\email{subhankarsau18@gmail.com}



\date{\today}
\dedicatory{}

\abstract In this paper, we classify the characteristic matrices associated to quasitoric manifolds over a vertex cut of a finite product of simplices.
We discuss the integral cohomology rings of these quasitoric manifolds with possibly minimal generators and show several relations among the product of these generators. 
We classify integral cohomology rings (up to isomorphism as graded rings) of the quasitoric manifolds over the vertex cut of a finite product of simplices.
\endabstract

\maketitle

\section{Introduction}

Davis and Januszkiewicz introduced a class of even-dimensional smooth manifolds and called them `toric manifolds' in their pioneering paper \cite{DJ}. These manifolds are the topological generalizations of smooth projective toric varieties, also known as toric manifolds, see \cite{M08}.
Buchstaber and Panov \cite{BP} used the term `quasitoric manifolds' instead of `toric manifolds' to avoid any confusion. Choi-Masuda-Suh \cite{CMS08} and Masuda-Panov \cite{MP08} also used the second naming `quasitoric manifold' for this class of manifolds. So we are choosing the second terminology without hesitation. Briefly, 
a $2n$-dimensional manifold $M^{2n}$ with an effective `locally standard' smooth action of the compact $n$-torus $T^n$ having the orbit space a simple $n$-polytope $P$ is called a quasitoric manifold.

Classification of toric manifolds is an active topic in mathematics nowadays; see \cite{MS08}, \cite{CMS08}, \cite{CMS10}, \cite{CM12} and \cite{Choi15}.
Choi-Masuda-Suh \cite{CMS08} classified the quasitoric manifolds over a finite product of simplices. In \cite{HKMP}, Hasui-Kuwata-Masuda-Park discussed and classified toric manifolds over a vertex cut of an $n$-dimensional cube.
Inspired by the above works, in this paper, we study some properties of quasitoric manifold whose orbit spaces are the product of $m$ many simplices with a vertex cut. 
In particular, if the simplices are closed intervals then their product is an $m$-cube. Consequently, a few of our results generalizes some results of \cite{HKMP}.

The paper is organized as follows. We start Section \ref{Sec_Quasitoric manifold} with the concept of simple polytopes $P$ and their codimension one faces $\mathcal{F}(P)$. Then we briefly recall the constructive definition of a quasitoric manifold, denoted by $X(P, \lambda)$, corresponding to a simple polytope $P$ and a characteristic function $\lambda \colon \mathcal{F}(P) \to \ZZ^{\dim P}$, following \cite{DJ} and \cite{BP}. 
The characteristic function induces a linear map $\ZZ^{|\mathcal{F}(P)|} \to \ZZ^{\dim P}$. The matrix $\Lambda$ of this linear map is called a characteristic matrix. It is enough to characterize the matrix $\Lambda$ up to an equivalence to classify quasitoric manifolds over $P$ up to an equivariant homeomorphism, see \cite[Proposition 1.8]{DJ} and \cite[Proposition 5.14]{BP}.
We study the quasitoric manifold over a finite product of simplices following \cite{CMS08}. For a vector matrix, we discuss the concept of principal minors and determinants. We also discuss the classification of square vector matrices $A$ with all the proper principal minors $1$ and their determinants belonging to $\{1, -1\}$. 
We define the distance, denoted by $D(u_1, u_0)$, between two vertices $u_0, u_1$ of a polytope. If $P$ is a product of $m$ many simplices and $\lambda' \colon \mathcal{F}(P) \to \ZZ^{\dim P}$ then we uniquely define the matrix $A_{\textbf{v}}$ at a vertex $\textbf{v}$ of $P$. 

In Section \ref{Sec_Quasitoric manifolds over vertex cut}, we first recall the idea of a vertex cut of a simple polytope. We denote the vertex cut of the product of simplices $P=\prod_{j=1}^m \Delta^{n_j}$ along a vertex $\tilde{\textbf{v}}$ by $\widebar{P}$ where $\Delta^{n_j}$ is simplex of dimension $n_j$ and $\tilde{\textbf{v}}$ is an $m$-distant vertex from a fixed point $\textbf{v}_0 \in V(P)$. 
If $\mathcal{F}(P)=\{F_1, \dots, F_r\}$ then $\mathcal{F}(\widebar{P})=\{F_1 \cap \widebar{P}, \dots, F_m \cap \widebar{P}, \widebar{F}\}$ where the codimension one face $\widebar{F}$ arises due to the vertex cut of $P$.
Let $$\widebar{\lambda} \colon \mathcal{F}(\widebar{P}) \to \ZZ^{\dim \widebar{P}}$$ be a characteristic function defined as in \eqref{Eq_characteristic function over widebarP}. 
One can define a matrix $A_{\textbf{u}}$ uniquely for each vertex $\textbf{u}$ of $\widebar{P}$ similar to the construction of $A_{\textbf{v}}$. 
The characteristic function $\widebar{\lambda}$ on $\widebar{P}$ induces a function $\lambda' \colon \mathcal{F}(P) \to \ZZ^{\dim P}$, see \eqref{Eq_induced lambda on P}. The function $\lambda'$ may not be a characteristic function on $P$; see Figure \ref{Fig_Example of induced lambda may not be a char fn.} for an example. Then, we prove one of the main results of our paper.

\begin{theorem}[Theorem \ref{Theorem on the matrix A}]
 Let $X(\widebar{P}, \widebar{\lambda})$ be a quasitoric manifold where $\widebar{P}$ is vertex cut at $\tilde{\textbf{\emph{v}}}$ of the polytope $P=\prod_{j=1}^m \Delta^{n_j} $ and $\widebar{\lambda}$ is defined as in \eqref{Eq_characteristic function over widebarP} satisfying 
 \begin{align}\label{Eq_next use 2}
     \det A_{\textbf{\emph{u}}}=\begin{cases}
-1 \quad \emph{if } D(\textbf{\emph{u}},\textbf{\emph{u}}_0)=\emph{odd}\\
+1 \quad \emph{if } D(\textbf{\emph{u}},\textbf{\emph{u}}_0)=\emph{even}
\end{cases}
 \end{align}
 for $\textbf{\emph{u}} \in V(\widebar{P})$ and $\textbf{\emph{u}}_0$ is a fixed vertex of $\widebar{P}$ corresponding to the vertex $\textbf{\emph{v}}_0$ of $P$.
Then the matrix $A_{\tilde{\textbf{\emph{v}}}}$ can be characterized based on its determinant.
\end{theorem}
We also show how the vector $\textbf{b} : = \overline{\lambda}(\overline{F})$, assigned to the facet $\overline{F}$ of $\widebar{P}$, is dependent on the vectors assigned to the other facets of $\widebar{P}$. Moreover, we classify the vector $\textbf{b}$ depending on the determinant of $A_{\tilde{\textbf{v}}}$, see Theorem \ref{Theorem on the vector b}.

In Section \ref{Sec_Cohomology of quasitoric manifolds over the vertex cut of a finite product of simplices}, we calculate the integral cohomology of toric manifold $X(\widebar{P},\widebar{\lambda})$ with possibly minimal generators where $\widebar{P}$ is the vertex cut along $\tilde{\textbf{v}}$ of the product of simplices $P = \prod_{j=1}^m \Delta^{n_j}$ and $\widebar{\lambda}$ is defined as in \eqref{Eq_characteristic function over widebarP} satisfying \eqref{Eq_next use 2}.
 We show several relations among the product of these generators. Let $ \overline{\lambda} \colon \mathcal{F}(\overline{P}) \to \ZZ^{\dim P}$ be a characteristic function and $\lambda \colon \mathcal{F}(P) \to \ZZ^{\dim P}$ be the function induced from $\overline{\lambda}$ as in \eqref{Eq_induced lambda on P}. Then, for $\det A_{\tilde{{\textbf{v}}}}=0$, we prove the following.
 
\begin{theorem}[Theorem \ref{Th_classify rings when det =0}]
Let $P=\prod_{j=1}^m \Delta^{n_j}$ be a product of simplices as in \eqref{Eq_P is a product of simplices} and $\widebar{P}$ a vertex cut of $P$ along the vertex $\tilde{\emph{\textbf{v}}}$ such that $\det A_{\tilde{\emph{\textbf{v}}}}=0$.
Then the cohomology rings $H^*(X(\widebar{P}, \widebar{\lambda}))$ are isomorphic to each other if $b_i=0$ for $i \neq N_s$ and $s=1, \dots, m$ in the vector $\textbf{b}$ where $N_s=\sum_{j=1}^s n_j$.
\end{theorem}
\noindent We provide an additive basis of $H^4(X(\widebar{P},\widebar{\lambda}))$, see Theorem \ref{Theorem for additive basis}. Also, for $\det A_{\tilde{{\textbf{v}}}}=(-1)^m$, we prove the following.
\begin{theorem}[Theorem \ref{Th_classify ring when det=1}]
Let $P=\prod_{j=1}^m \Delta^{n_j}$ be a product of simplices as in \eqref{Eq_P is a product of simplices} with $m \geq 2$, $n \geq 3$ and $\widebar{P}$ a vertex cut of $P$ along the vertex $\tilde{\emph{\textbf{v}}}$ such that $\det A_{\tilde{\emph{\textbf{v}}}}=(-1)^m$.
Then $H^*(X(\widebar{P}, \widebar{\lambda}))$ and $H^*(X(\widebar{P}, \widebar{\lambda}'))$ are isomorphic as graded rings if and only if $H^*(X(P, \lambda))$ and $H^*(X(P, \lambda'))$ are isomorphic as graded rings where $\lambda$ and $\lambda'$ are characteristic functions induced from $\widebar{\lambda}$ and $\widebar{\lambda}'$ respectively.
\end{theorem}

\section{Quasitoric manifold over a finite product of simplices}\label{Sec_Quasitoric manifold}

In this section, we recall the constructive definition and some properties of polytopes and quasitoric manifolds following \cite{BP}. 
We also recall the classification of a square matrix $A$ which has all the proper principal minors $1$. Then we define the distance between two vertices of a polytope. Given a function $\lambda \colon \mathcal{F}(P) \to \ZZ^{\dim P}$ on a finite product of simplices $P$, we define a matrix $A_{\textbf{v}}$ at a vertex $\textbf{v} \in V(P)$ of the product of simplices $P$ whose columns are the vectors assigned to the facets intersecting at $\textbf{v}$ with an ordering fixed uniquely.

An $n$-dimensional simple polytope is an $n$-dimensional convex polytope such that at each vertex, exactly $n$ codimension one faces (called facets) intersect. Let $Q$ be an $n$-dimensional simple polytope in $\RR^n$. 
We denote the facets and the vertices of $Q$ by $\mathcal{F}(Q)$ and $V(Q)$ respectively. Let $\mathcal{F}(Q):=\{F_1, \dots, F_r\}$.

\begin{definition}\label{linear independency of char vec}
Let $\lambda \colon \mathcal{F}(Q) \to \mathbb{Z}^n$ be a function such that 
\begin{multline}\label{Eq_char vec unimodular}
\{\lambda(F_{i_1}),\dots,\lambda(F_{i_k})\} \text{ is a part of } \ZZ \text{-basis of } \ZZ^n \text{ whenever } F_{i_1} \cap \dots \cap F_{i_k}\neq \varnothing.
\end{multline} 
Then $\lambda$ is called a {characteristic function} on $Q$. We denote $\lambda(F_i)$ by $\lambda_i$, the {characteristic vector} assigned to the facet $F_i$ for $i=1, \dots,r$. The pair $(Q,\lambda)$ is called a {characteristic pair}. Also, the $(n \times r)$ matrix $\Lambda=\big( \lambda_1 ~~ \lambda_2 ~~ \dots~~ \lambda_r \big)$ is called characteristic matrix associated to the characteristic function $\lambda$ where $\lambda_i$ is the $i$-th column of $\Lambda$ for $i=1, \dots, r$.
\end{definition}


We briefly recall the construction of a quasitoric manifold from a characteristic pair $(Q,\lambda)$. Let $p$ be a point in $Q$ belongs to the relative interior of a codimension $d$ face $F$ of $Q$. Thus $F=\bigcap_{j=1}^d F_{i_j}$ for some unique facets $F_{i_1}, \dots, F_{i_d}$ of $Q$. We consider the standard $n$-dimensional torus $T^n:=(\RR \otimes_{\ZZ} \ZZ^n)/\ZZ^n$ and $\mathbb{R}^n$ as the Lie algebra of $T^n$. Note that each $\lambda_i$ determines a unique line in $\mathbb{R}^n$. We denote the image of this line under the map $exp \colon \mathbb{R}^n \to T^n$ by $T_i$ for $i=1, \dots,r$.
Let $T_F$ be the $d$-dimensional subtorus of $T^n$ generated by $\{T_{i_1}, \ldots, T_{i_d}\}$.

Consider the identification $\sim$ on $T^n \times Q$ defined as follows:
\begin{equation*}\label{eq_equivalence_rel}
(t,p) \sim (s,q)~ \text{if and only if}~ p=q ~\text{and}~ t^{-1}s \in T_F,
\end{equation*}
where $p$ is a point belonging to the relative interior of the face $F$. Then the quotient space 
\begin{equation*}
X(Q,\lambda)  :=(T^n \times Q) / \sim
\end{equation*}
is a $T^n$-manifold and is called a quasitoric manifold, see \cite[Chapter 5]{BP}. The $T^n$ action is induced by the group operation on the first factor of $T^n \times Q$ and the orbit map 
\begin{equation}\label{Eq_define pi}
\pi \colon X(Q, \lambda) \to Q
\end{equation}
is given by $\pi((t,p)_{\sim}) = p$ where $(t,p)_{\sim}$ is the equivalence class of $(t,p)$. 
We note that if $M$ is a quasitoric manifold then there is a characteristic pair $(Q, \lambda)$ such that $M$ is equivariantly homeomorphic to $X(Q, \lambda)$, see \cite[Proposition 1.8]{DJ}.


Next, we discuss some properties of the characteristic functions on a finite product of simplices following \cite{CMS08}. Let 
\begin{equation}\label{Eq_P is a product of simplices}
P :=\prod_{j=1}^m \Delta^{n_j} 
\end{equation}
be a product of $m$ simplices of dimensions $n_1, \dots, n_{m-1} $and $n_m$. So, the dimension of $P$ is $\sum_{j=1}^m n_j$ which we denote by $n$. 
Let us denote 
\begin{equation}\label{Eq_define Ns}
N_s:= \sum_{j=1}^s n_j
\end{equation}
for $s=1, \dots, m$. Note that $N_1=n_1$ and $N_m=n$. Also, we assume $N_0:=0$.

Let $V(\Delta^{n_j}):=\{v_0^j, \dots ,v_{n_j}^j\}$ be the vertex set of $\Delta^{n_j}$ for $j=1, \dots,m$. So, 
\begin{equation}\label{Eq_vertex set of P}
V(P)=\{v_{k_1 k_2 \dots k_m}:=(v_{k_1}^1, v_{k_2}^2, \dots, v_{k_m}^m )~ |~0\leq k_j \leq n_j, j=1, \dots, m\}.
\end{equation} 
For each $j \in \{1, \dots, m\}$, let $\mathcal{F}(\Delta^{n_j}):=\{f_0^j, \dots f_{n_j}^j\}$ where the facet $f_{k_j}^j$ does not contain the vertex $v_{k_j}^j$ in $\Delta^{n_j}$ for $0 \leq k_j \leq n_j$. Then the facets of $P$ are given by
\begin{align}\label{Eq_facet set of P}
\mathcal{F}(P) & :=\{F_{k_j}^j:=\Delta^{n_1} \times \dots \times \Delta^{n_{j-1}} \times f_{k_j}^j \times \Delta^{n_{j+1}}\times \dots \times \Delta^{n_m} ~|~0\leq k_j \leq n_j, 1 \leq j \leq m\}.
\end{align}
Observe that the simple polytope $P$ has $\prod_{j=1}^m n_j$ many vertices and $(n+m)$ many facets. Note that the vertex $v_{k_1 k_2 \dots k_m}$ in $P$ is the unique intersection of the $n$-many facets $\mathcal{F}(P) \setminus \{F_{k_j}^j ~|~j=1, \dots,m\}$. In particular,
\begin{equation}\label{at v00000}
v_{0 \dots 0}= F_1^1 \cap \dots \cap F_{n_1}^1 \cap \dots \cap F_1^m \cap \dots \cap F_{n_m}^m.
\end{equation}

We call two faces $F_1, F_2 \not \in V(P)$  adjacent if $F_1 \cap F_2 \neq \varnothing$. Note that the facets $F_1^1 , \dots , F_{n_1}^1 , \dots , F_1^m, \dots , F_{n_m}^m$ are adjacent to the vertex $v_{0\dots0}$. Let
\begin{equation}\label{Eq_define lambda on P}
\lambda \colon \mathcal{F}(P) \to \mathbb{Z}^n.
\end{equation}
be a characteristic function on $P$. 
Then $\{\lambda(F^1_1), \dots, \lambda(F^1_{n_1}), \dots,\lambda(F^m_1), \dots, \lambda(F^m_{n_m})\}$ is a $\ZZ$-basis of $\ZZ^n$.
So, we may assume that it is the standard basis vectors; that is \begin{align}\label{Eq_lambda at origin}
\lambda(F_1^1)=e_1, &\dots, \lambda(F^1_{n_1})=e_{n_1},\\ &\vdots \nonumber \\
\lambda(F_1^j)=e_{N_{j-1}+1}, & \dots, \lambda(F_{n_j}^j)=e_{N_j},\nonumber \\ & \vdots \nonumber
\\ \lambda(F^m_1)=e_{N_{m-1}+1}, & \dots, \lambda(F^m_{n_m})=e_n.\nonumber
\end{align}
For the remaining $m$ facets $F^1_0, \dots, F^m_0$ of $P$, we assign
 \begin{equation} \label{Eq_lambda to other facets}
\textbf{a}_j:=\lambda(F_0^j) \in \mathbb{Z}^n \quad ~\text{for}~j=1, \dots, m.
 \end{equation}
This gives the vector matrices of order $1 \times m$ and $m \times m$, which we can also write as an $n \times m$ scalar matrix as follows.
 \begin{align}\label{Eq_vector matrix}
 A:=
\begin{pmatrix} 
 \textbf{a}_1 & \textbf{a}_2 & \dots & \textbf{a}_m
\end{pmatrix}_{1 \times m} =
\begin{pmatrix}
\textbf{a}_1^1 & \dots & \textbf{a}_m^1 \\
\vdots & \dots & \vdots \\
\textbf{a}_1^m & \dots & \textbf{a}_m^m 
\end{pmatrix}_{m \times m}=
\begin{pmatrix}
a_{11}^1 & \dots & a_{m1}^1\\
\vdots & \dots & \vdots \\
a^1_{1n_1} & \dots & a^1_{mn_1}\\
\vdots & \dots & \vdots \\
a^m_{11} & \dots & a^m_{m1} \\
\vdots & \dots & \vdots \\
 a^m_{1n_m} & \dots & a^m_{mn_m} 
\end{pmatrix}_{n \times m}
\end{align} 
where $\textbf{a}_j \in \mathbb{Z}^n$ is considered as the $j$-th column vector of $A$, $\textbf{a}_j^k \in \mathbb{Z}^{n_k}$ is the $(k,j)$-th entry of the $m \times m$ vector matrix and $a_{ji}^k \in \ZZ$ is the $(N_{k-1}+i,j)$-th entry of the $n \times m$ scalar matrix.
For simplicity of notation, we often denote the $\ell$-th coordinate of the vector $\textbf{a}_j$ by $a_{j \ell}$ where $\ell=N_{k-1}+i$ for some $1 \leq i \leq n_{k}$ and $1 \leq k \leq m$.

 \begin{definition}
 Given an $n \times n$ scalar matrix $M$, a $k$-th($1 \leq k < n$) proper principal minor of $M$ is the determinant of the $k \times k$ submatrix obtained by removing $(n-k)$ many rows and the corresponding columns from $M$. 
 \end{definition}



Let $A$ be an $n \times m$ scalar matrix as in \eqref{Eq_vector matrix}. Given $ k_j \in \{1 , \dots, n_j\}$ for $j=1, \dots, m$, we define an $m \times m$ submatrix 
\begin{equation}\label{Eq_submatrix}
A_{k_1, \dots, k_m}:=
\begin{pmatrix}
a^1_{1k_1} & \dots  & a^1_{mk_1} \\
\vdots & \ddots & \vdots\\
a^m_{1k_m}  & \dots & a^m_{mk_m}
\end{pmatrix}_{m \times m}
\end{equation}
of $A$ whose $j$-th row is the $(N_{j-1}+k_j)$-th row of $A$ where $N_{j-1}$ is defined in \eqref{Eq_define Ns}. 
By proper principal minor of the matrix $A$ we mean the proper principal minors of the matrix $A_{k_1, \dots, k_m}$ for all $ k_j \in \{1 , \dots, n_j\}$ and $j=1, \dots, m$.

\begin{proposition}\cite[Lemma 3.2]{CMS08}\label{l 3.2 cms08}
Let $\lambda$ be a characteristic function on $P=\prod_{j=1}^m \Delta^{n_j}$ as in \eqref{Eq_define lambda on P} and $A$ the matrix associated to $\lambda$ as in \eqref{Eq_vector matrix}.
Then the condition on $\lambda$ in \eqref{Eq_char vec unimodular} is equivalent to all principal minors of $A$ is $\pm 1$.
\end{proposition}


\begin{remark}
Proposition \ref{l 3.2 cms08} implies that the diagonal entry vector $$\textbf{a}^j_j=(a^j_{j1}, \dots, a^j_{jn_j}) \in \ZZ^{n_j}$$ of $A$ is equipped with entries $a^j_{ji}=\pm 1$ for $i=1, \dots, n_j$ and $j=1, \dots, m$. 
\end{remark}

Recall that an upper or lower triangular matrix is called \emph{unipotent} if all of its diagonal entries are units.
The $i$-th row (respectively column) of the vector matrix $A_{m \times m}$ defined in \eqref{Eq_vector matrix} is called \emph{elementary row} (respectively \emph{elementary column}) if $\textbf{a}^j_j=\textbf{1}=(1,\dots,1)$ (respectively $\textbf{1}^t$) and all other entries of the row (respectively column) are $\textbf{0}=(0, \dots, 0)$ vectors.
\begin{definition}\label{def_dets of A}
Let $A$ be a matrix of the form \eqref{Eq_vector matrix}. By `all the determinants of $A$', we denote the determinants of all the $m \times m$ submatrices $A_{k_1, \dots, k_m}$ of $A$ for some $N_{j-1} < k_j \leq N_j$ and $j=1, \dots, m$ where $A_{k_1, \dots, k_m}$ is defined in \eqref{Eq_submatrix}. 
\end{definition}

\begin{proposition}\label{lemma on matrix}
Let $A$ be an $m \times m$ vector matrix of the form \eqref{Eq_vector matrix} such that all the proper principal minors of $A$ are $1$. If all the determinants of $A$ are $1$, then $A$ is conjugate to a unipotent upper triangular vector matrix of the following form:
\begin{equation}\label{Eq_matrix for 1st case}
\begin{pmatrix}
\textbf{1} & \textbf{a}_2^1 & \textbf{a}_3^1 & \dots & \textbf{a}_m^1 \\
\textbf{0} & \textbf{1} & \textbf{a}_3^2 & \dots & \textbf{a}_m^2 \\
\vdots & \ddots & \ddots & \ddots & \vdots \\
\textbf{0} & \dots & \dots & \textbf{1} & \textbf{a}_{m}^{m-1} \\
\textbf{0} & \dots & \dots & \textbf{0} & \textbf{1}
\end{pmatrix}
\end{equation}
where $\textbf{0}=(0, \dots, 0)^t$, $\textbf{1}=(1, \dots, 1)^t$ of appropriate size. If all the determinants of $A$ are $\pm 1$ and at least one of them is $-1$, then $A$ is conjugate to a vector matrix of the following form:
\begin{equation}\label{Eq_matrix for 2nd case}
\begin{pmatrix}
\textbf{1} & \textbf{a}^1_2 & \textbf{0} & \dots & \textbf{0} \\
\textbf{0} & \textbf{1} & \textbf{a}^2_3 & \dots & \textbf{0} \\
\vdots & \vdots & \ddots & \ddots & \vdots \\
\textbf{0} & \textbf{0} & \dots & \textbf{1} & \textbf{a}^{m-1}_m \\
\textbf{a}^m_1 & \textbf{0} & \dots & \textbf{0} & \textbf{1}
\end{pmatrix}
\end{equation}
where $\textbf{a}_1^m$ and $\textbf{a}_{j}^{j-1}$ are non-zero vectors for $j=2, \dots, m$ such that $$\prod_{j=2}^m a^{j-1}_{ji} \times a^m_{1i}=(-1)^m 2,$$ $a^{j-1}_{ji}$ is any non-zero coordinate of $\textbf{a}_{j}^{j-1}$ for $j=2, \dots, m$ and $a^m_{1i}$ is any non-zero coordinate of $\textbf{a}_1^m$. Moreover, the non-zero coordinates in $\textbf{a}_1^m, \textbf{a}_2^1, \dots, \textbf{a}_m^{m-1}$ are all same and they are either $\pm 1$ or $\pm 2$.
\end{proposition}

\begin{proof}
The proof of the proposition follows from Lemma 3.3 of \cite{MP08} and Lemma 5.1 of \cite{CMS08}.
\end{proof}



The edges of a polytope $Q$ form a connected graph and one can have the following.

\begin{definition}[Distance function]\label{Def_distance fn}
Let $Q$ be a polytope and $v_1$, $v_2$ two different vertices in $Q$. A \emph{path} between $v_1$ and $v_2$ is a sequence of edges $\xi_1, \xi_2, \dots, \xi_d$ such that $v_1 \in \xi_1$, $v_2 \in \xi_d$ and $\xi_i \cap \xi_{i+1}$ is a vertex of both for $i=1, \dots, (d-1)$.
The \emph{distance} between two vertices $v_1$ and $v_2$ is the minimum $d$ and it is denoted by $D(v_1,v_2)$.
If $D(v_1,v_2)=d \geq 1$, then we call $v_2$ a $d$-distant vertex from $v_1$ in $Q$. 
\end{definition}

We remark that if $P$ is a product of $m$ simplices as in \eqref{Eq_P is a product of simplices} and $v_1, v_2 \in V(P)$, then $D(v_1, v_2) \leq m$. Also, there may exist more than one different path of length $m$ if $D(v_1, v_2) =m$.

\begin{figure}
\begin{tikzpicture}[scale=.4]
\draw[dashed] (0,0)--(1,2)--(3,2)--(4,0)--(2,-2)--cycle;
\draw (0,5)--(1,7)--(3,7)--(4,5)--(2,3)--cycle;
\draw [dashed] (1,2)--(1,7);
\draw[dashed] (3,2)--(3,7);
\draw (4,0)--(4,5);
\draw (2,-2)--(2,3);
\draw (0,0)--(0,5);
\draw (0,0)--(2,-2)--(4,0);
\node [right] at (4,5) {$v_1$};
\node [right] at (4,0) {$v_2$};
\node [left] at (0,0) {$v_3$};

\node at (12,5) {$D(v_1,v_2)=1$};
\node at (12,4) {$D(v_1,v_3)=3$};
\node at (12,3) {$D(v_3,v_2)=2$};
\node at (12.7,2) {$D(v_i,v_i)=0~\forall~i$};
\end{tikzpicture}
\caption{A distance function on the vertices of a polytope.}\label{Fig_Distance function}
\end{figure}
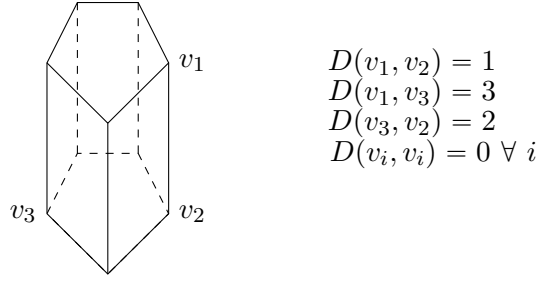

Let $(P, \lambda)$ be a characteristic pair where $P$ is a product of $m$ simplices as in \eqref{Eq_P is a product of simplices} and $\lambda$ a characteristic function on $P$ as defined in \eqref{Eq_define lambda on P} satisfying \eqref{Eq_lambda at origin} and \eqref{Eq_lambda to other facets}.
Let $\textbf{v}$ be a vertex in $P$. So, $\textbf{v}=\bigcap_{j=1}^n F_j$ for some unique facets $F_j$'s of $P$. 
Without any loss of generality, at $\textbf{v}_0:=v_{0 \dots 0} \in V(P)$, we may fix the order of columns as follows.
\begin{align}\label{Eq_A_v0}
A_{\textbf{v}_0}:&=\begin{pmatrix} \lambda(F_1^1) & \dots & \lambda(F^1_{n_1}) &  \dots & \lambda(F_1^m) & \dots  & \lambda(F^m_{n_m}) \end{pmatrix}
\\
&=\begin{pmatrix} e_1 & \dots & e_{n_1} & \dots & e_{N_{m-1}+1} & \dots & e_n \end{pmatrix} \nonumber
\end{align}
where $e_1, \dots, e_n$ are the standard $\ZZ$-basis of $\ZZ^n$. 

Now for any vertex $\textbf{v} \in V(P)$, we fix the ordering of the column entries in $A_{\textbf{v}}$ as follows. Let $D(\textbf{v},\textbf{v}_0)=d >0$. By \eqref{Eq_vertex set of P}, $\textbf{v}=v_{\ell_1 \ell_2 \dots \ell_m}$ for some $0 \leq \ell_j \leq n_j$, $j=1, \dots, m$. 
Then $$\textbf{v}=\bigcap_{\substack{j=1\\ k_j \neq \ell_j}}^m F^j_{k_j}$$ where $F^j_{k_j}$'s are defined in \eqref{Eq_facet set of P}. 
Note that if $\ell_j \neq 0$ for $j\in \{1, \dots, m\}$, then $e_{N_{j-1}+\ell_j}$ is replaced by $\textbf{a}_j$ by keeping the order of other columns of $A_{\textbf{v}_0}$ intact where $N_{j-1}$ is defined in \eqref{Eq_define Ns}.

Let us consider a path of length $d$ from $\textbf{v}_0$ to $\textbf{v}$. Notice that at each step along the path from $\textbf{v}_0$ to $\textbf{v}$ exactly one vector is replaced and the new vector takes the place of the vector which it replaces.
That is if $\xi_1, \dots, \xi_d$ is the sequence of edges joining $\textbf{v}_0$ to $\textbf{v}$ such that $\textbf{v}_0 \in \xi_1$, $\textbf{v} \in \xi_d$ and $\xi_i \cap \xi_{i+1}=\textbf{v}_i$ then the matrix $A_{\textbf{v}_{i+1}}$ is formed by replacing exactly one standard basis vector in the columns of $A_{\textbf{v}_{i}}$ for $i=1, \dots, d-1$.
Note that the matrix $A_{\textbf{v}}$ does not alter by the choice of the path if we choose any other shortest path of length $d$.
If $\textbf{v}$ is a vertex such that $D(\textbf{v}, \textbf{v}_0)=m$, i.e., $\ell_j \neq 0$ for all $j=1, \dots, m$.
Then the matrix $A_{\textbf{v}}$ is given by
\begin{multline*}
A_{\textbf{v}}=\big( e_1~ ~ \dots ~~ e_{\ell_1-1} ~~ \textbf{a}_1 ~~ e_{\ell_1+1} ~~ \dots ~~ e_{N_1} ~~ e_{N_1+1} ~~ \dots ~~ e_{N_1+\ell_2-1} ~~ \textbf{a}_2 ~~ e_{N_1+\ell_2+1} ~~ \dots \\ e_{N_2} ~~ \dots~ ~e_{N_{m-1}+1} ~~ \dots ~~ e_{N_{m-1}+\ell_{m}-1} ~~ \textbf{a}_m ~~ e_{N_{m-1}+\ell_{m}+1} ~~ \dots ~~ e_{N_m} \big).
\end{multline*}
Here we want to notify that throughout this paper, the vectors $e_i$ and $\textbf{a}_j$ of $\ZZ^n$ are considered as the columns of matrices $A_{\textbf{v}}$'s for some $i=1, \dots, n$ and $j=1, \dots, m$.
In the following, we illustrate this process for a particular case.

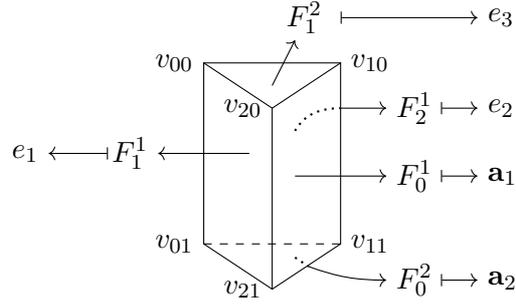
\begin{figure}
\begin{tikzpicture}[scale=.6]
\draw (0,0)--(1.5,-1)--(3,0)--(3,4)--(1.5,3)--(0,4)--cycle;
\draw (1.5,-1)--(1.5,3);
\draw (0,4)--(3,4);
\draw [dashed] (0,0)--(3,0);

\node [left] at (0,4) {$v_{00}$};
\node [right] at (3,4) {$v_{10}$};
\node [left] at (1.5,3) {$v_{20}$};
\node [left] at (0,0) {$v_{01}$};
\node [right] at (3,0) {$v_{11}$};
\node [left] at (1.5,-1) {$v_{21}$};

\draw [->] (2,1.5)--(4,1.5);
\node [right] at (4,1.5) {$F_0^1$};
\draw [|->] (5.2,1.5)--(6,1.5);
\node [right] at (6,1.5) {$\textbf{a}_1$};

\draw[dotted, thick] (2,2.5) to [in=180, out=60] (3,3);
\draw[->] (3,3)--(4,3);
\node [right] at (4,3) {$F_2^1$};
\draw [|->] (5.2,3)--(6,3);
\node [right] at (6,3) {$e_2$};

\draw[->] (1,2)--(-1,2);
\node [left] at (-1,2) {$F_1^1$};
\draw [|->] (-2.1,2)--(-3.4,2);
\node [left] at (-3.4,2) {$e_1$};

\draw [dotted, thick] (2,-.3)--(2.3,-.5);
\draw [->] (2.3,-.5) to [out=340, in=180] (4,-.8);
\node [right] at (4,-.8) {$F_0^2$};
\draw [|->] (5.2,-.8)--(6,-.8);
\node [right] at (6,-.8) {$\textbf{a}_2$};

\draw [->] (1.5,3.5)--(2,4.5);
\node at (2.2,5) {$F_1^2$};
\draw [|->] (3,5)--(6,5);
\node [right] at (6,5) {$e_3$};

\end{tikzpicture}
\caption{The characteristic function $\lambda$ on a prism following \eqref{Eq_define lambda on P}.} \label{Fig_Figure to understand how the vectors are replaced at vertices.}
\end{figure}

\begin{example}\label{Eg_matrices at vertices}
Consider the $3$-dimensional prism $P = \Delta^2 \times I$ where $\Delta^2$ is $2$-simplex and $I$ is the interval as in Figure \ref{Fig_Figure to understand how the vectors are replaced at vertices.}.
Then
\begin{equation*}
V(P)=\{v_{00},v_{01},v_{10},v_{11},v_{20},v_{21}\}, \quad \mbox{and} \quad \mathcal{F}(P)=\{F_0^1, F_1^1, F_2^1, F_0^2, F_1^2\}.
\end{equation*}
Define a characteristic function $\lambda \colon \mathcal{F}(P) \to \ZZ^3$ following \eqref{Eq_define lambda on P}. Thus
$$ \lambda(F_1^1)=e_1,\quad \lambda(F_2^1)=e_2, \quad \lambda(F_1^2)=e_3,\quad \lambda(F_0^1)=\textbf{a}_1,\quad \lambda(F_0^2)=\textbf{a}_2.$$
Then the matrix $A_{v_{00}}=(e_1,e_2,e_3)$. Now consider the vertex $v_{01}$, then $D(v_{01},v_{00})=1$. So, the matrix $A_{v_{01}}=(e_1,e_2,\textbf{a}_2)$ is formed by replacing one vector of $A_{v_{00}}$. 
Similarly, we get $A_{v_{10}}=(e_1, \textbf{a}_1, e_3),$ and $A_{v_{20}}=(\textbf{a}_1,e_2, e_3)$.
The vertices $\{v_{11}, v_{21}\}$ are at distance $2$ from $v_{00}$. 
So $A_{v_{11}}$ is formed by replacing two vectors of $A_{v_{00}}$. The consistency of ordering of the vectors in $A_{v_{11}}$ can be followed by using $D(v_{11},v_{01})=1$. Therefore, $A_{v_{11}}=(e_1, \textbf{a}_1,\textbf{a}_2)$. Similarly, we get $A_{v_{21}}=(\textbf{a}_1, e_2, \textbf{a}_2)$. \qed
\end{example}

\section{Quasitoric manifolds over a vertex cut of a finite product of simplices} \label{Sec_Quasitoric manifolds over vertex cut}

The idea of blowup of a polytope along a face is studied in \cite{GP} and \cite{BSSau21}. The vertex cut of a polytope is a blowup of a polytope along a vertex. In this section, first we recall the definition of blowup and vertex cut. 
Then we characterise quasitoric manifolds, or equivalently characteristic pairs, over a vertex cut of a product of finitely many simplices.

\begin{definition}\label{def_blowup polytope}
Let $Q$ be an $n$-dimensional simple polytope in $\mathbb{R}^n$ and $F$ a proper face of $Q$. Take a hyperplane $H$ in $\mathbb{R}^n$ such that one open half space $H_{>0}$ contains all the vertices of $V(F)$ and $V(Q) \setminus V(F)$ is a subset of the other open half space $H_{<0}$. Then the polytope $\widebar{Q}:=Q \cap H_{\leq 0}$ is called a blowup of $Q$ along the face $F$.
\end{definition}

If $F=\textbf{v} \in V(Q)$ in Definition \ref{def_blowup polytope}, then the blowup is called a \emph{vertex cut} of $Q$ along the vertex $\textbf{v}$.
The new facet $\widebar{F}$ in $\widebar{Q}$, which arises due to vertex cut of $Q$ is called the facet associated to $\textbf{v}$. Then $\widebar{F}$ is an $(n-1)$-dimensional simplex.

Let $P$ be the product $\prod_{j=1}^m \Delta^{n_j}$ as in \eqref{Eq_P is a product of simplices}. The vertices and facets of $P$ are
\begin{align*}
V(P)&=\{v_{k_1 k_2 \dots k_m}:=(v_{k_1}^1, v_{k_2}^2, \dots, v_{k_m}^m )~ |~0\leq k_j \leq n_j, j=1, \dots m\} \quad \text{and}\\
\mathcal{F}(P)&=\{F_{k_j}^j ~|~0\leq k_j \leq n_j, 1 \leq j \leq m\},
\end{align*}
as discussed in \eqref{Eq_vertex set of P} and \eqref{Eq_facet set of P} respectively. Let $\widebar{P}$ be a vertex cut of $P$ at the vertex $\tilde{\textbf{v}}:=v_{n_1n_2 \dots n_m}$. Note that the distance $D(\textbf{v}_0,\tilde{\textbf{v}})=m$ where $\textbf{v}_{0}=v_{0 \dots 0}$ in $P$.
Then the vertex set and the facet set of $\widebar{P}$ are respectively
\begin{align}\label{Eq_vertex and facet set of bar P}
    V(\widebar{P})&:=(V(P) \setminus \{\tilde{\textbf{v}}\} ) \cup V(\widebar{F}) \quad \text{and}\\
    \mathcal{F}(\widebar{P})&:=\{\widebar{F}^j_{k_j} := F^j_{k_j} \cap \widebar{P}~ |~ F^j_{k_j} \in \mathcal{F}(P)\} \cup \{\widebar{F}\}. \nonumber
\end{align}
For convenience of notation, when we consider $v_{k_1 k_2 \dots k_m}$ as a vertex of $ \widebar{P}$, we denote it by $u_{k_1 k_2 \dots k_m}$ if $(k_1, k_2, \dots, k_m) \neq (n_1, n_2, \dots, n_m)$.
We denote the vertex $u_{00 \dots 0}$ by $\textbf{u}_0$. 

 Let
 \begin{equation}\label{Eq_characteristic function over widebarP}
     \widebar{\lambda} \colon \mathcal{F}(\widebar{P}) \to \ZZ^n
 \end{equation}
be a characteristic function on $\overline{P}$. Note that $\textbf{u}_0=\widebar{F}_1^1 \cap \dots \cap \widebar{F}_{n_1}^1 \cap \dots \cap \widebar{F}_1^m \cap \dots \cap \widebar{F}_{n_m}^m$. So, one can assume that
 \begin{align*}
  \widebar{\lambda}(\widebar{F}^1_j)&:=e_j \quad \quad \text{for } j=1, \dots,n_1,\\ \vdots \nonumber \\
  \widebar{\lambda}(\widebar{F}^m_j)&:=e_{N_{m-1}+j} \quad \text{for } j=1, \dots,n_m, \nonumber
 \end{align*}
 Let 
  \begin{equation}
 \widebar{\lambda}(\widebar{F}) :=\textbf{b} \in \ZZ^n ~~\mbox{and }~~  \widebar{\lambda}(\widebar{F}^j_0) := \textbf{a}_j \in \ZZ^n \quad \text{for } j=1, \dots,m.
  \end{equation}
One can define a matrix $A_{\textbf{u}}$ uniquely for each vertex $\textbf{u}$ of $\widebar{P}$ similar to the construction of $A_{\textbf{v}}$ for each vertex $\textbf{v}$ of $P$. 
Then the matrix associated to the vertex $\textbf{u}_0$ is $$A_{\textbf{u}_0}=\begin{pmatrix} e_1 & \dots & e_{N_1} & e_{N_1+1}& \dots & e_{N_{m-1}} & \dots & e_{N_m}\end{pmatrix}$$ where $e_i$ is considered as the $i$-th column of the matrix $A_{\textbf{u}_0}$.
Recall from Section \ref{Sec_Quasitoric manifold} that the $\ell$-th coordinate of the vector $\textbf{a}_j$ is denoted by $a_{j \ell}$ where $\ell=N_{k-1}+i$ for some $1 \leq i \leq n_{k}$ and $k=1, \dots, m$.

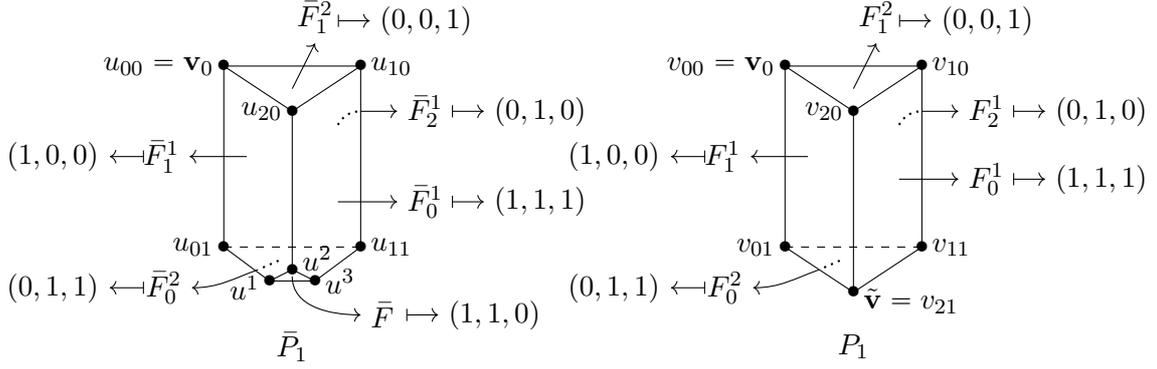
\begin{figure}
\begin{center}
\begin{tikzpicture}[scale=.6]
\draw (0,0)--(0,4)--(3,4)--(3,0)--(2,-.75)--(1,-.75)--(1.5,-.5)--(1.5,3)--(3,4);
\draw (0,0)--(1,-.75);
\draw (0,4)--(1.5,3);
\draw[dashed] (0,0)--(3,0);
\draw (1.5,-.5)--(2,-.75);
\node at (1.5,-2.2) {$\widebar{P}_1$};

\draw[dotted, thick] (2.5,2.7) to [out=60, in=180] (3,3);
\draw [->] (3,3)--(3.8,3);
\node[right] at (3.8,3) {$\widebar{F}_2^1$};
\draw[|->] (5,3)--(5.7,3);
\node[right] at (5.7,3) {$(0,1,0)$};

\draw [->] (.5,2)--(-.75,2);
\node [left] at (-.75,2) {$\widebar{F}_1^1$};
\draw[|->] (-1.75,2)--(-2.5,2);
\node[left] at (-2.5,2) {$(1,0,0)$};

\draw [->] (2.5,1)--(3.8,1);
\node[right] at (3.8,1) {$\widebar{F}_0^1$};
\draw[|->] (5,1)--(5.7,1);
\node[right] at (5.7,1) {$(1,1,1)$};

\draw [->] (1.5,3.5)--(2,4.5);
\node[above] at (2,4.5) {$\widebar{F}_1^2$};
\draw[|->] (2.5,5)--(3.2,5);
\node[right] at (3.2,5) {$(0,0,1)$};

\draw [dotted, thick] (1.2,-.3)--(.75,-.5);
\draw [->] (.75,-.5) to [out=200, in=360] (-.7,-.9);
\node[left] at (-.7,-.9) {$\widebar{F}_0^2$};
\draw[|->] (-1.75,-.9)--(-2.5,-.9);
\node[left] at (-2.5,-.9) {$(0,1,1)$};

\draw[->] (1.5,-.65) to [out=270, in=180] (3,-1.5);
\node[right] at (3,-1.5) {$\widebar{F}$};
\draw[|->] (4,-1.5)--(4.7,-1.5);
\node[right] at (4.7,-1.5) {$(1,1,0)$};

\node at (0,0) {$\bullet$};
\node[left] at (0,0) {$u_{01}$};

\node at (0,4) {$\bullet$};
\node[left] at (0,4) {$u_{00}=\textbf{v}_0$};

\node at (1.5,3) {$\bullet$};
\node[left] at (1.5,3) {$u_{20}$};

\node at (3,4) {$\bullet$};
\node[right] at (3,4) {$u_{10}$};

\node at (3,0) {$\bullet$};
\node[right] at (3,0) {$u_{11}$};

\node at (1,-.75) {$\bullet$};
\node[left] at (1,-.9) {$u^1$};

\node at (1.5,-.5) {$\bullet$};
\node[right] at (1.5,-.25) {$u^2$};

\node at (2,-.75) {$\bullet$};
\node[right] at (2,-.75) {$u^3$};

\begin{scope}[xshift=350]
\draw (0,0)--(0,4)--(3,4)--(3,0)--(1.5,-1)--(1.5,3)--(3,4);
\draw (0,0)--(1.5,-1);
\draw (0,4)--(1.5,3);
\draw[dashed] (0,0)--(3,0);

\draw[dotted, thick] (2.5,2.7) to [out=60, in=180] (3,3);
\draw [->] (3,3)--(3.8,3);
\node[right] at (3.8,3) {$F^1_2$};
\draw[|->] (5,3)--(5.7,3);
\node[right] at (5.7,3) {$(0,1,0)$};

\draw [->] (.5,2)--(-.75,2);
\node [left] at (-.75,2) {$F^1_1$};
\draw[|->] (-1.75,2)--(-2.5,2);
\node[left] at (-2.5,2) {$(1,0,0)$};

\draw [->] (2.5,1.5)--(3.8,1.5);
\node[right] at (3.8,1.5) {$F^1_0$};
\draw[|->] (5,1.5)--(5.7,1.5);
\node[right] at (5.7,1.5) {$(1,1,1)$};

\draw [->] (1.5,3.5)--(2,4.5);
\node[above] at (2,4.5) {$F_1^2$};
\draw[|->] (2.5,5)--(3.2,5);
\node[right] at (3.2,5) {$(0,0,1)$};

\draw [dotted, thick] (1.2,-.3)--(.75,-.5);
\draw [->] (.75,-.5) to [out=200, in=360] (-.7,-.9);
\node[left] at (-.7,-.9) {$F_0^2$};
 \draw[|->] (-1.75,-.9)--(-2.5,-.9);
\node[left] at (-2.5,-.9) {$(0,1,1)$};

\node at (0,0) {$\bullet$};
\node[left] at (0,0) {$v_{01}$};

\node at (0,4) {$\bullet$};
\node[left] at (0,4) {$v_{00}=\textbf{v}_0$};

\node at (1.5,3) {$\bullet$};
\node[left] at (1.5,3) {$v_{20}$};

\node at (3,4) {$\bullet$};
\node[right] at (3,4) {$v_{10}$};

\node at (3,0) {$\bullet$};
\node[right] at (3,0) {$v_{11}$};

\node at (1.5,-1) {$\bullet$};
\node[right] at (1.5,-1.2) {$\tilde{\textbf{v}}=v_{21}$};
\node at (1.5,-2.2) {$P_1$};
\end{scope}
\end{tikzpicture}
\end{center}
\caption{A vertex cut of a prism.}
\label{Fig_Example of induced lambda may not be a char fn.}
\end{figure}

 The characteristic pair $(\widebar{P}, \widebar{\lambda})$ induces a map $\lambda \colon \mathcal{F}(P) \to \ZZ^n$ defined by 
 \begin{equation}\label{Eq_induced lambda on P}
     \lambda(F^j_{k_j}):=\widebar{\lambda}(\widebar{F}_{k_j}^j) 
 \end{equation}
for $j=1, \dots, m$ and $1 \leq k_j \leq n_j$. In general, this map $\lambda$ may not be a characteristic function on $P$. For example, in Figure \ref{Fig_Example of induced lambda may not be a char fn.}, $\widebar{P}_1$ is the vertex cut of the prism $P_1$ along $\tilde{\textbf{v}}=F_1^1 \cap F^1_0 \cap F^2_0$ and $\{(1,0,0), (1,1,1), (0,1,1)\}$ is not a $\ZZ$-basis of $\ZZ^3$. 


Consider the vertex $\tilde{\textbf{v}}:=v_{n_1 \dots n_m} \in V(P)$. Then from Section \ref{Sec_Quasitoric manifold}, \begin{equation}\label{Eq_tilde v}
\tilde{\textbf{v}}=F^1_1 \cap \dots \cap F^1_{n_1-1} \cap F^1_0 \cap F^2_1 \cap \dots \cap F_1^m \cap \cdots \cap F^m_{n_m-1} \cap F^m_0.
\end{equation}
So we can associate the following matrix 
  \begin{align}\label{Eq_matrix A_tilde v}
     A_{\tilde{\textbf{v}}}:
     =(e_1  ~~ \dots ~~ e_{N_1-1} ~~ \textbf{a}_1 ~~ e_{N_1+1} ~~ \dots ~~ e_{N_{(m-1)}-1} ~~ \textbf{a}_{m-1} ~~ e_{N_{(m-1)}+1} ~~ \dots ~~ e_{N_m-1} ~~ \textbf{a}_m)
  \end{align}
to the vertex $\tilde{\textbf{v}}$ with the ordering of the columns according to the ordering of the facets in the intersection in \eqref{Eq_tilde v}. 
 Let $A'_{\tilde{\textbf{v}}}$ be the submatrix of $A_{\tilde{\textbf{v}}}$ defined as in \eqref{Eq_submatrix} for $k_j=n_j$ and $j=1, \dots, m$. So, 
 \begin{equation}\label{Eq_matrix A' tilde v}
 A'_{\tilde{\textbf{v}}}=\begin{pmatrix}
 a_{1N_1} & a_{2N_1} & a_{3N_1} & \dots & a_{mN_1}\\
 a_{1N_2} & a_{2N_2} & a_{3N_2} & \dots & a_{mN_2}\\
 \vdots & \vdots & \vdots & \dots & \vdots \\
 a_{1N_m} & a_{2N_m} & a_{3N_m} & \dots & a_{mN_m}
 \end{pmatrix}_{m \times m}
 \end{equation}
where $N_1, N_2, \dots, N_m$ are defined in \eqref{Eq_define Ns}. Also, from \eqref{Eq_matrix A_tilde v}, we have $$\det A_{\tilde{\textbf{v}}}= \det A'_{\tilde{\textbf{v}}}.$$ 
 The following theorem characterizes the matrix $A_{\tilde{\textbf{v}}}$ corresponding to the characteristic function $\widebar{\lambda}$ on $\widebar{P}$. 

\begin{theorem}\label{Theorem on the matrix A}
 Let $X(\widebar{P}, \widebar{\lambda})$ be a quasitoric manifold where $\widebar{P}$ is vertex cut at $\tilde{\textbf{\emph{v}}}=v_{n_1 \dots n_m}$ of the polytope $P=\prod_{j=1}^m \Delta^{n_j} $ and $\widebar{\lambda}$ is defined as in \eqref{Eq_characteristic function over widebarP} satisfying 
 \begin{align}\label{Eq_determinant at vertices}
     \det A_{\textbf{\emph{u}}}=\begin{cases}
-1 \quad \emph{if } D(\textbf{\emph{u}},\textbf{\emph{u}}_0)=\emph{odd}\\
+1 \quad \emph{if } D(\textbf{\emph{u}},\textbf{\emph{u}}_0)=\emph{even}
\end{cases}
 \end{align}
 for $\textbf{\emph{u}} \in V(\widebar{P})$.
Then the matrix $A_{\tilde{\textbf{\emph{v}}}}$ can be characterized based on its determinant.
\end{theorem}

\begin{proof}
The $(n \times m)$ matrix $A=\begin{pmatrix} \textbf{a}_1 & \textbf{a}_2 & \dots & \textbf{a}_m \end{pmatrix}$ is a submatrix of the $(n \times n)$ matrix $A_{\tilde{\textbf{v}}}$. We can consider $A$ as an $(m \times m)$ vector matrix as discussed in \eqref{Eq_vector matrix}. Since the columns of $A_{\tilde{\textbf{v}}}$ which are not columns of $A$ are from the standard basis vectors of $\ZZ^n$, it is enough to understand $A$ to characterize $A_{\tilde{\textbf{v}}}$.

We compute the determinants of $A_{{u}_{k_1 k_2 \dots k_m}}$ for $(k_1, k_2, \dots, k_m) \neq (0, 0, \dots, 0),$ and $ (k_1, k_2, \dots, k_m) \neq (n_1, n_2, \dots, n_m)$ where ${u}_{k_1 k_2 \dots k_m} \in V(\widebar{P})$. Recall that $0 \leq k_j \leq n_j$ for $j=1, \dots, m$.
Since $P$ is a finite product of simplices, $D({u}_{k_1 k_2 \dots k_m}, \textbf{u}_0)=d$ for $ 0< d \leq m$ if and only if $k_j \neq 0$ for exactly $d$ many $j$'s. 
This implies 
\begin{equation}\label{Eq_det A_u}
\det A_{u_{k_1 k_2 \dots k_m}}=(-1)^d
\end{equation}
from our hypothesis \eqref{Eq_determinant at vertices}. Note that if $k_j \neq 0$, then the standard basis vector $e_{N_{j-1}+k_j}$ in $A_{\textbf{u}_0}=I_n$ is replaced by $\textbf{a}_j$ to form $A_{u_{k_1 k_2 \dots k_m}}$ where $I_n$ is the $n \times n$ identity matrix.
In particular, if we consider $k_1 \in \{1, \dots, n_1\}$ and $k_j=0$ for $j=2, \dots, m$, then the corresponding vertex is $u_{0 \dots 010 \dots0} \in V(\widebar{P})$. So $D(\textbf{u}_0, u_{0 \dots 010 \dots0})=1$ where 1 is at the $k_1$-th position. Then $$A_{{u}_{0 \dots 010 \dots0}}:=(e_1, \dots, e_{k_1-1,}\textbf{a}_1, e_{k_1+1}, \dots, e_n).$$
The hypothesis in \eqref{Eq_determinant at vertices} implies $\det A_{{u}_{0 \dots 010 \dots0}}=-1$. Thus $a_{1k_1}=-1$ where $a_{1k_1}$ is the $k_1$-th entry of $\textbf{a}_1$. Using a similar argument for all one distant vertices, the hypothesis in \eqref{Eq_determinant at vertices} leads us to
\begin{equation}\label{Eq_a^j_jk}
a_{j \ell}=-1
\end{equation}
for $N_{j-1} < \ell \leq N_j$ and $j=1, \dots,m$ where $a_{j \ell}$ is the $\ell$-th coordinate of the vector $\textbf{a}_j$.

 The $d$-th proper principal minors of $A$ are the determinants of the matrices $A_{u_{k_1 k_2 \dots k_m}}$ where exactly $d$ many entries of $(k_1, k_2, \dots, k_m)$ are non-zero. Thus from \eqref{Eq_det A_u}, we have the $d$-th proper principal minors of $A$ are $(-1)^d$. Let us consider $$B=-A.$$ 
 Then $B$ is an $(m \times m)$ vector matrix with $d$-th proper principal minors $(-1)^d \times d$-th proper principal minors of $A$. Thus all the proper principal minors of the matrix $B$ are $1$. 
Now we classify the matrix $B$ using Proposition \ref{lemma on matrix} up to a conjugation. 

\textbf{Case 1:} Let $\det A_{\tilde{\textbf{v}}}=(-1)^m$. 
Also, when all the entries of $(k_1, k_2, \dots, k_m)$ are non-zero and $(k_1, k_2, \dots, k_m) \neq (n_1, n_2, \dots, n_m)$ then $\det A_{u_{k_1 k_2 \dots k_m}}=(-1)^m$, follows from \eqref{Eq_det A_u}.
Thus the determinants of the vector matrix $A$ are $(-1)^m$, where the determinants of $A$ are defined in Definition \ref{def_dets of A}.  
Let $B_m$ be an $m \times m$ submatrix of $B$.  So $B_m=-A_{u_{k_1k_2 \dots k_m}}$ for some $(k_1, k_2, \dots, k_m)$ with all the entries non-zero and 
\begin{align*}
\text{det} B_m = (-1)^m \times  \text{det}A_{u_{k_1k_2 \dots k_m}}=(-1)^{2m}=1.
\end{align*}
Thus $B=-A$ is an upper triangular matrix of the form \eqref{Eq_matrix for 1st case} and
\begin{equation}\label{Eq_matrix when det is 1}
    A_{\tilde{\textbf{v}}}=
    \left( \begin{array}{*{20}c}
        1 & 0 & \dots & 0 & -1 & 0 & \dots & a_{21} & \dots & a_{31} & \dots & 0 & a_{m1}\\
        0 & 1 & \dots & 0 & -1 & 0 & \dots & a_{22} & \dots & a_{32} & \dots & 0 & a_{m2}\\
        \vdots & \vdots & \ddots & \vdots & \vdots & \vdots & \ddots & \vdots & \ddots & \vdots & \ddots & \vdots & \vdots\\
        0 & 0 & \dots & 1 & -1 & 0 & \dots & a_{2N_1-1} & \dots & a_{3N_1-1} & \dots & 0 & a_{mN_1-1}\\
        0 & 0 & \dots & 0 & -1 & 0 & \dots & a_{2N_1} & \dots & a_{3N_1} & \dots & 0 & a_{mN_1}\\
        0 & 0 & \dots & 0 & 0 & 1 & \dots & a_{2N_1+1} & \dots & a_{3N_1+1} & \dots &  0 &a_{mN_1+1}\\
        \vdots & \vdots & \ddots & \vdots & \vdots & \vdots & \ddots & \vdots & \ddots & \vdots & \ddots & \vdots & \vdots\\
        \vdots & \vdots & \ddots & \vdots & \vdots & \vdots & \ddots & \vdots & \ddots & \vdots & \ddots & \vdots & \vdots\\
        0 & 0 & \dots & 0 & 0 & 0 & \dots & 0 & \dots & 0 & \dots & 1 & -1\\
        0 & 0 & \dots & 0 & 0 & 0 & \dots & 0 & \dots & 0 & \dots & 0 & -1\\
    \end{array} \right).
\end{equation}

\textbf{Case 2:} Let $\det A_{\tilde{\textbf{v}}} \neq (-1)^m$. This implies at least one of the determinants of $B$ is not $1$. Thus by Proposition \ref{lemma on matrix}, $A=-B$ is of the followinh form:
\begin{equation}\label{Eq_classify A case 2}
A=\begin{pmatrix}
\textbf{-1} & \textbf{a}^1_2 & \textbf{0} & \dots & \textbf{0} \\
\textbf{0} & \textbf{-1} & \textbf{a}^2_3 & \dots & \textbf{0} \\
\vdots & \vdots & \ddots & \ddots & \vdots \\
\textbf{0} & \textbf{0} & \dots & \textbf{-1} & \textbf{a}^{m-1}_m \\
\textbf{a}^m_1 & \textbf{0} & \dots & \textbf{0} & \textbf{-1}
\end{pmatrix}.
\end{equation}
Then the matrix $A'_{\tilde{\textbf{v}}}$ defined in \eqref{Eq_matrix A' tilde v} can be given by
\begin{equation*}
 A'_{\tilde{\textbf{v}}}=\begin{pmatrix}
 -1 & a_{2N_1} & 0 & \dots & 0 & 0\\
 0 & -1 & a_{3N_2} & \dots & 0 & 0\\
 \vdots & \vdots & \vdots & \dots & \vdots &\vdots \\
 0 & 0 & 0 & \dots & -1 & a_{mN_{m-1}}\\
 a_{1N_m} & 0 & 0 & \dots & 0 & -1
 \end{pmatrix}_{m \times m}
 \end{equation*}
and
\begin{equation}\label{Eq_determinant of Atildev}
    \det A_{\tilde{\textbf{v}}}=\det A'_{\tilde{\textbf{v}}}=(-1)^m + (-1)^{m-1}a_{1N_m}a_{2N_1} \dots a_{mN_{m-1}}.
\end{equation}
Since $\det A_{\tilde{\textbf{v}}} \neq (-1)^m$, then $a_{jN_{j-1}} \neq 0$ for $j=2, \dots, m$ and $a_{1N_m} \neq 0$.

Now we consider the matrices $A_{u_{k_1 k_2 \dots k_m}}$ with all the entries of $(k_1, k_2, \dots, k_m)$ are non-zero and $(k_1, k_2, \dots, k_m) \neq (n_1, n_2, \dots, n_m)$. Note that these matrices have $A$ as their submatrix whereas the other columns are appropriate standard basis vectors of $\ZZ^n$. Then using \eqref{Eq_classify A case 2} and hypothesis \eqref{Eq_determinant at vertices}, we have
$$\text{det} A_{u_{k_1 k_2 \dots k_m}}=(-1)^m+(-1)^{m-1} a_{1\ell_1} a_{2 \ell_2} \dots a_{m \ell_m}=(-1)^m\Rightarrow a_{1\ell_1} a_{2 \ell_2} \dots a_{m \ell_m}=0$$
where $\ell_1=N_{m-1}+k_1$ and $\ell_j=N_{j-1}+k_j$ for $j=1, \dots, (m-1)$. 
Let us consider $k_j=n_j$ for $j=1, \dots, i-1, i+1, \dots, m$ and $k_i>0$. So, the vertex $u_{n_1 \dots n_{i-1}k_i n_{i+1} \dots n_m}$ is at a distance $m$ from $\textbf{u}_0$ in $\widebar{P}$. 
Then $A_{u_{n_1 \dots n_{i-1}k_i n_{i+1} \dots n_m}}=(-1)^m$. Thus $$a_{1N_m}a_{2N_1} \dots a_{(i-1)N_{i-2}} a_{i(N_{i-1}+k_i)} a_{(i+1)N_i} \dots a_{mN_{m-1}}=0.$$
Since the product $a_{1N_m}a_{2N_1} \dots a_{(i-1)N_{i-2}}  a_{(i+1)N_i} \dots a_{mN_{m-1}}$ is non-zero from \eqref{Eq_determinant of Atildev}, we have $a_{i(N_{i-1}+k_i)}=0$.
Similar arguments for the vertices ${u_{k_1 k_2 \dots k_m}}$ with all the entries of $(k_1, k_2, \dots, k_m)$ non-zero and $(k_1, k_2, \dots, k_m) \neq (n_1, n_2, \dots, n_m)$,  we have $$a_{j(N_{j-1}+k_j)}=0$$ for $0< k_j < n_j$ and $j=1, \dots, m$.
Thus, we get
\begin{equation}\label{Eq_matrix when det is not 1}
    A_{\tilde{\textbf{v}}}=
    \left( \begin{array}{*{20}c}
        1 & 0 & \dots & 0 & -1 & 0 & \dots & 0 & \dots  & 0 & 0\\
        0 & 1 & \dots & 0 & -1 & 0 & \dots & 0 & \dots  & 0 & 0\\
        \vdots & \vdots & \ddots & \vdots & \vdots & \vdots & \ddots & \vdots & \ddots  & \vdots & \vdots\\
        0 & 0 & \dots & 1 & -1 & 0 & \dots & 0 & \dots  & 0 & 0 \\
        0 & 0 & \dots & 0 & -1 & 0 & \dots & a_{2N_1} & \dots  & 0 & 0\\
        0 & 0 & \dots & 0 & 0 & 1 & \dots & -1& \dots  &  0 & 0 \\
        \vdots & \vdots & \ddots & \vdots & \vdots & \vdots & \ddots & \vdots & \ddots  & \vdots & \vdots\\
         0 & 0 & \dots & 0 & 0 & 0 & \dots & 0 & \dots  & 0 & a_{m N_{m-1}}\\
        \vdots & \vdots & \ddots & \vdots & \vdots & \vdots & \ddots & \vdots & \ddots  & \vdots & \vdots\\
        0 & 0 & \dots & 0 & 0 & 0 & \dots & 0 & \dots  & 1 & -1\\
        0 & 0 & \dots & 0 & a_{1N_m} & 0 & \dots & 0 & \dots  & 0 & -1\\
    \end{array} \right).
\end{equation}

\end{proof}
 
 The following result states how the vector $\textbf{b}~(=\widebar{\lambda}(\widebar{F}))$ is dependent on the columns of $A_{\tilde{\textbf{v}}}$ where $\widebar{F}$ is the facet of $\widebar{P}$ associated to the vertex $\tilde{\textbf{v}} \in P$.

 \begin{theorem}\label{Theorem on the vector b}
 Let $X(\widebar{P}, \widebar{\lambda})$ be a quasitoric manifold where $\widebar{P}$ is the vertex cut at $\tilde{\textbf{\emph{v}}}={v}_{n_1 \dots n_m}$ of the polytope $P=\prod_{j=1}^m \Delta^{n_j} $ and $\widebar{\lambda}$ is defined as in \eqref{Eq_characteristic function over widebarP} satisfying 
 \begin{align}\label{Eq_determinant at u}
     \det A_{\textbf{\emph{u}}}=\begin{cases}
-1 \quad \emph{if } D(\textbf{\emph{u}},\textbf{\emph{u}}_0)=\emph{odd}\\
+1 \quad \emph{if } D(\textbf{\emph{u}},\textbf{\emph{u}}_0)=\emph{even}
\end{cases}
 \end{align}
 for any $\textbf{\emph{u}} \in V(\widebar{P})$ and $\textbf{\emph{u}}_0:={u}_{0,\dots,0}$.
 Then we can determine \emph{$\textbf{b}$} according to the determinant of the matrix $A_{\tilde{\textbf{\emph{v}}}}$ as follows.
 
\textbf{Case 1:} If $\det A_{\tilde{\textbf{\emph{v}}}}=0$ then $$\sum_{j=1}^m b_{N_j}=-1$$ where $N_1,\dots, N_m$ are defined in \eqref{Eq_define Ns} and $b_i$ can be arbitrary if $i \not \in \{N_1, \dots, N_m\}$.

\textbf{Case 2:} If $\det A_{\tilde{\textbf{\emph{v}}}} \neq 0$ then $$b_i = \frac{(-1)^m}{\det A_{\tilde{\textbf{\emph{v}}}}} \sum_{q=1}^n A_{(i,q)}$$  for $i=1, \dots, n$ where $A_{(i,q)}$ is the $(i,q)$-th entry of the matrix $A_{\tilde{\textbf{\emph{v}}}}$. 
 \end{theorem}
 
 \begin{proof}
Let $V(\widebar{F})=\{u^1, \dots , u^n\}$. Also, $u^j=E_j \cap \widebar{F}$ for unique edge $E_j$ of $P$ adjacent to $\tilde{\textbf{v}}$.
There are $m$ many shortest paths from $\textbf{v}_0$ to $\tilde{\textbf{v}}$ in $P$. 
A vertex on the last edge of the shortest path gives an $m$-distant vertex in $\widebar{F} \subset \widebar{P}$.
Continuing similarly, we get $m$ many $m$-distant vertices in $\{u^1, \dots , u^n\}$. 
Without loss of generality, we may assume that the $m$-distant vertices are $\{u^1, \dots, u^m\}$ and the other $(n-m)$ vertices $\{u^{m+1}, \dots, u^n\}$ are at a distance $m+1$. Thus using \eqref{Eq_determinant at u}, we get
 \begin{align}
     \det A_{\textbf{u}}=
     \begin{cases}
     (-1)^m \quad &\text{if } \textbf{u}\in \{u^1, \dots, u^m\}, \\
     (-1)^{m+1} \quad &\text{if } \textbf{u}\in \{u^{m+1}, \dots, u^n\}.
     \end{cases}
 \end{align}
 
 Now we focus on the matrices $A_{u^1}, \dots, A_{u^m}$ and observe their formation. First, consider the vertex $u^1$ and the matrix $A_{u^1}$. 
From Section \ref{Sec_Quasitoric manifold}, we recall $$\tilde{\textbf{v}}=v_{n_1 \dots n_m}= \bigcap_{j=1}^m ( \bigcap_{k_j=0}^{n_j-1} F^j_{k_j} ).$$ Also the matrix at $A_{\tilde{\textbf{v}}}$ at the vertex $\tilde{\textbf{v}}$ is defined in \eqref{Eq_matrix A_tilde v} following the ordering of the columns.
Let $$S_{\tilde{\textbf{v}}}:=\{ F^j_{k_j} \colon 1 \leq j \leq m, ~0 \leq k_j \leq (n_j-1) \}$$ be the set of facets containing $\tilde{\textbf{v}}$. Then $$E_{1}= \bigcap_{F_i \in S_{\tilde{\textbf{v}}}  \setminus \{F^1_0\}} F_i$$ in $P$. Thus in $\widebar{P}$, we have
 $$u^1= \bigcap_{F_i \in S_{\tilde{\textbf{v}}}  \setminus \{F^1_0\}} \widebar{F}_i \cap \widebar{F}.$$ 
Then the matrix $A_{u^1}$ for the vertex $u^1$ is defined by replacing $\textbf{a}_1$ by $\textbf{b}$; that is
 \begin{equation}\label{Eq_matrix at u^1}
 A_{u^1}=\left( \begin{array}{*{20}c}
 e_1 & \dots & e_{n_1-1} & \textbf{b} & e_{n_1+1} & \dots & \textbf{a}_{m-1} & e_{N_{m-1}+1} & \dots & e_{n-1}  & \textbf{a}_m
  \end{array} \right).
 \end{equation}
 
Similarly, the matrix $A_{u^j}$ is formed by replacing the column vector $\textbf{a}_j$ by $\textbf{b}$ in the matrix $A_{\tilde{\textbf{v}}}$; that is 
 \begin{multline}\label{Eq_matrix at vertices of bar P}
 A_{u^j}=\big(
 e_1 ~~ \dots ~~ e_{N_1-1} ~~ \textbf{a}_1 ~~ e_{N_1+1} ~~ \dots ~~ e_{N_j-1} ~~ \textbf{b} ~~ e_{N_{(j+1)}+1} ~~ \dots ~~ e_{N_m-1} ~~ \textbf{a}_m
 \big)
 \end{multline}
 for $j=2, \dots, m$. So, $\det A_{u^j}=(-1)^m$ for $j=1, \dots, m$.

Now $D(u^p, u^q)=1$ for $u^p, u^q \in V(\widebar{F})$ and $p \neq q$. 
So $$D(u^j, u^{m+\alpha})=1$$ for $j=1, \dots, m$ and $\alpha=1, \dots, (n-m)$. Thus the matrix $A_{u^{m + \alpha}}$ associated to the vertex $u^{m + \alpha}$ is obtained by replacing the column $e_{N_{j-1} + \alpha}$ by $\textbf{a}_j$ in $A_{u^j}$ defined in \eqref{Eq_matrix at vertices of bar P} where $N_{j-1} < \alpha < N_j$ and $j=1, \dots, m$; that is
\begin{multline}\label{Eq_matrix at m+1 distant vertices of bar P}
 A_{u^{m+\alpha}}=\big( e_1 ~~ \dots ~~ e_{N_1-1} ~~ \textbf{a}_1 ~~ e_{N_1+1} ~~ \dots ~~ \textbf{a}_{j-1} ~~ \dots e_{N_{j-1}+\alpha-1} ~~\textbf{a}_{j} ~~  e_{N_{j-1}+\alpha+1} ~~\\  \dots ~~ e_{N_j-1} ~~ \textbf{b} ~~ e_{N_{(j+1)}+1} ~~ \dots ~~ e_{N_m-1} ~~ \textbf{a}_m \big).
 \end{multline} 
Since $D(\textbf{u}_0, u^{m+\alpha})=m+1$, we have $\det A_{u^{m+\alpha}}= (-1)^{m+1}$ for $\alpha=1, \dots, (n-m)$. Interchanging the columns $\textbf{a}_1$ and $\textbf{b}$ in $A_{u^{m+1}}$ leads to 
\begin{multline}\label{Eq_column interchange}
\det \big( e_1 ~~ \dots ~~ e_{N_1-1} ~~ \textbf{a}_1 ~~ e_{N_1+1} ~~ \dots ~~ \textbf{a}_{j-1} ~~ \dots e_{N_{j-1}+\alpha-1} ~~ \textbf{b} ~~ e_{N_{j-1}+\alpha+1}\\  \dots ~~ e_{N_j-1} ~~ \textbf{a}_{j} ~~ e_{N_{(j+1)}+1}~~ \dots~~ e_{N_m-1}~~ \textbf{a}_m \big) = (-1)^{m+2}=(-1)^m.
\end{multline}

Let $\widetilde{A}_{\tilde{\textbf{v}}}$ be the adjoint matrix of $A_{\tilde{\textbf{v}}}$ and $\widetilde{A}_{(p,q)}$ the $(p,q)$-th entry of $\widetilde{A}_{\tilde{\textbf{v}}}$. Then 
\begin{equation}\label{Eq_relation between cofactor}
\widetilde{A}_{(p,q)}=(-1)^{(p+q)} \times \text{cofactor of } a_{(q,p)}
\end{equation}
where $a_{(q,p)}$ denotes the $(q,p)$-th entry of the matrix $A_{\tilde{\textbf{v}}}$.
From \eqref{Eq_matrix at vertices of bar P}, \eqref{Eq_matrix at m+1 distant vertices of bar P} and \eqref{Eq_column interchange}, we get $$\sum_{p=1}^n \widetilde{A}_{(p,q)} b_p=(-1)^m$$ for $q=1, \dots, n$.
Thus 
\begin{equation}\label{Eq_A tilde b}
    \widetilde{A}_{\tilde{\textbf{v}}}\textbf{b}=(-1)^m \textbf{1}
\end{equation}
where $\textbf{1}$ is the matrix $(1, \dots, 1)_{n \times 1}^t$ and considering $\textbf{b}$ as column vector. Also note that $$A_{\tilde{\textbf{v}}} \widetilde{A}_{\tilde{\textbf{v}}}=(\det A_{\tilde{\textbf{v}}}) I_n$$ where $I_n$ is the $n \times n$ identity matrix.

\textbf{Case1:} Let $\det A_{\tilde{\textbf{v}}}=0$. Then $A_{\tilde{\textbf{v}}} \widetilde{A}_{\tilde{\textbf{v}}}=\textbf{0}_{n \times n}$ where $\textbf{0}_{n \times n}$ is an $n \times n$ zero matrix.
From Theorem \ref{Theorem on the matrix A}, the matrix $A_{\tilde{\textbf{v}}}$ is of the form \eqref{Eq_matrix when det is not 1}.
We first show that in $\widetilde{A}_{\tilde{\textbf{v}}}$, the columns other than $N_1, N_2, \dots, N_m$-th columns have only zeros as its entries. 
First, we concentrate on the $s$-th column of $\widetilde{A}_{\tilde{\textbf{v}}}$ where $s \neq N_j$ for any $j=1, \dots, m$. Its entries are calculated following \eqref{Eq_relation between cofactor} when $q=s$. 
As $s \neq N_j$ for $j=1, \dots, m$ then $e_s$ must be a column vector in $A_{\tilde{\textbf{v}}}$. The $(s,s)$-th entry of $A_{\tilde{\textbf{v}}}$ is $1$ and its cofactor is $0$ as $\det A_{\tilde{\textbf{v}}}=0$. Thus $\widetilde{A}_{(s,s)}=0$.

Now we look at $\widetilde{A}_{(t,s)}$, that is the cofactor of $a_{(s,t)}$ in $A_{\tilde{\textbf{v}}}$ where $t \neq s$. 
Let $M_{(s,t)}$ be the cofactor matrix of $a_{(s,t)}$ in $A_{\tilde{\textbf{v}}}$ where $t \neq s$. If $t <s$, then $(s-1)$-th column of $M_{(s,t)}$ is a zero column. The $s$-th column of $M_{(s,t)}$ is zero column if $t>s$. This zero column comes from the zeros of $e_s$. Thus $$\widetilde{A}_{(t,s)}=0$$ when $s \neq N_j$ for $j=1, \dots, m$.
So, the columns other than $N_1, N_2, \dots, N_m$ have only zeros as their entries. 

Let us denote the $i$-th row vector of $\widetilde{A}_{\tilde{\textbf{v}}}$ by $\textbf{r}_i$. Then from $A_{\tilde{\textbf{v}}} \widetilde{A}_{\tilde{\textbf{v}}}=\textbf{0}_{n \times n}$ and \eqref{Eq_matrix when det is not 1} we have 
\begin{equation*}
\adjustbox{max width=\textwidth}{$
\left( \begin{array}{*{20}c}
        1 & 0 & \dots & 0 & -1 & 0 & \dots & 0 & \dots  & 0 & 0\\
        0 & 1 & \dots & 0 & -1 & 0 & \dots & 0 & \dots  & 0 & 0\\
        \vdots & \vdots & \ddots & \vdots & \vdots & \vdots & \ddots & \vdots & \ddots  & \vdots & \vdots\\
        0 & 0 & \dots & 1 & -1 & 0 & \dots & 0 & \dots  & 0 & 0 \\
        0 & 0 & \dots & 0 & -1 & 0 & \dots & a_{2N_1} & \dots  & 0 & 0\\
        0 & 0 & \dots & 0 & 0 & 1 & \dots & -1& \dots  &  0 & 0 \\
        \vdots & \vdots & \ddots & \vdots & \vdots & \vdots & \ddots & \vdots & \ddots  & \vdots & \vdots\\
         0 & 0 & \dots & 0 & 0 & 0 & \dots & 0 & \dots  & 0 & a_{m N_{m-1}}\\
        \vdots & \vdots & \ddots & \vdots & \vdots & \vdots & \ddots & \vdots & \ddots  & \vdots & \vdots\\
        0 & 0 & \dots & 0 & 0 & 0 & \dots & 0 & \dots  & 1 & -1\\
        0 & 0 & \dots & 0 & a_{1N_m} & 0 & \dots & 0 & \dots  & 0 & -1\\
    \end{array} \right)    
    \begin{pmatrix}
    \textbf{r}_1 \\\textbf{r}_{2}\\ \vdots \\ \textbf{r}_{N_1-1} \\ \textbf{r}_{N_1} \\ \textbf{r}_{N_1+1} \\ \vdots \\ \textbf{r}_{N_{m-1}}\\ \vdots \\ \textbf{r}_{N_m-1} \\ \textbf{r}_{N_m}
    \end{pmatrix}
    =\textbf{0}_{ n \times n}.
    $}
\end{equation*}
Thus we get
\begin{align*}
\textbf{r}_{i} - \textbf{r}_{N_{j}}&=0 \quad \text{ for } N_{j-1} < i < N_j \text{ and } j=1, \dots, m,\\
\textbf{r}_{N_j} - a_{j'N_j} \textbf{r}_{N_{j'}}&=0 \quad \text{ for }  j=1, \dots, m-1,\\
\textbf{r}_{N_m} - a_{1N_m} \textbf{r}_{N_1} &=0 ,
\end{align*}
where $j'=j+1$.
We can rewrite \eqref{Eq_A tilde b} by $(-1)^m \widetilde{A}_{\tilde{\textbf{v}}} \textbf{b}=\textbf{1}$.
 Now we go through the following row operations on the augmented matrix $[(-1)^m\widetilde{A}_{\tilde{v}} ~|~ \textbf{1}]$
\begin{align*}
\textbf{r}_{i} - \textbf{r}_{N_{j}} \quad &\text{ for } N_{j-1} < i < N_j \text{ and } j=1, \dots, m \\
\textbf{r}_{N_j} - a_{j'N_j} \textbf{r}_{N_{j'}} \quad &\text{ for }  j=1, \dots, m-1
\end{align*}
 where $j'=j+1$. Then we get
\begin{align*}
\left(\begin{array}{ccccc|c}  
0 & \dots & \dots & 0 &0 & 0\\
\vdots & \ddots & \ddots & \vdots & \vdots & \vdots\\
0 & \dots & \dots & 0 &0 & 0\\
0 & \dots & \dots & 0 &0 & 1-a_{2N_1}\\
0 & \dots & \dots & 0 &0 & 0\\
\vdots & \ddots & \ddots & \vdots & \vdots & \vdots\\
0 & \dots & \dots & 0 &0 & 0\\
0 & \dots & \dots &0 &0 & 1-a_{mN_{m-1}}\\
0 & \dots & \dots & 0 &0 & 0\\
\vdots & \ddots & \ddots & \vdots & \vdots & \vdots\\
0 & \dots & \dots & 0 &0 & 0\\
* & \dots & \dots & *& (-1)^{m+1} & 1
\end{array}\right)_{n \times (n+1)}
\end{align*}
where the $*$'s in the last row are some integers. So, $\text{rank}\widetilde{A}_{\tilde{\textbf{v}}}=1$.
There exists a vector $\textbf{b}$ such that $(-1)^m \widetilde{A}_{\tilde{\textbf{v}}} \textbf{b}=\textbf{1}$ iff $\text{rank} [(-1)^m \widetilde{A}_{\tilde{\textbf{v}}}, \textbf{1}]=\text{rank} \widetilde{A}_{\tilde{\textbf{v}}}$.
Therefore, 
\begin{equation}\label{Eq_coordinates are 1}
a_{jN_{j-1}}=1 \quad \text{for } j=2, \dots, m.
\end{equation}
As $\det {A}_{\tilde{\textbf{v}}}=(-1)^m + (-1)^{m-1}a_{1N_m}a_{2N_1} \dots a_{mN_{m-1}}=0$, we have $a_{1N_m}=1$.
Thus 
\begin{multline*}
\det \big( e_1 ~~ \dots ~~ e_{N_1-1} ~~ \textbf{a}_1 ~~e_{N_1+1} ~~ \dots ~~ e_{N_{m-1}-1} ~~ \textbf{a}_{{m-1}} ~~ e_{N_{m-1}+1} ~~ \dots ~~ e_{N_m-1} ~~ \textbf{b} \big) \\=(-1)^{m+1} \sum_{j=1}^m b_{N_j}.
\end{multline*}
Also left hand side is the determinant of a matrix associated to a vertex in $\widebar{F}$ at a distance $m$ from $\textbf{u}_0$ in $\widebar{P}$. Therefore, $(-1)^{m+1} \sum_{j=1}^m b_{N_j}=(-1)^m$, by \eqref{Eq_determinant at u}. Then we conclude 
\begin{equation}
\sum_{j=1}^m b_{N_j}=-1
\end{equation}
and $b_i$ can be arbitrary if $i \not \in \{N_1, \dots, N_m\}$.

\textbf{Case2:} Let $\det A_{\tilde{\textbf{v}}} \neq 0$. Then from \eqref{Eq_A tilde b} and the fact $A_{\tilde{\textbf{v}}} \widetilde{A}_{\tilde{\textbf{v}}} =(\det A_{\tilde{\textbf{v}}}) I_n$, we have 
\begin{equation}
b_i = \frac{(-1)^m}{\det A_{\tilde{\textbf{v}}}} \sum_{q=1}^n a_{(i,q)}
\end{equation}
where $a_{(i,q)}$ is the $(i,q)$-th entry of the matrix $A_{\tilde{\textbf{v}}}$ for $i\in \{1, \dots, n\}$.

\textbf{Subcase 1:} Let $\det A_{\tilde{\textbf{v}}}=(-1)^m$. First, let $i \neq N_j$ for $j=1, \dots, m$. If $i > N_{j-1}$ for some $j=1, \dots, m-1$, then $$b_i= \sum_{\substack{k>j}}^m a_{ki}$$ for $i=1, \dots, N_{m-1}-1$ and $b_i=0$ for $N_{m-1} <i <N_m$.
Now if $i = N_{j}$ for some $j=1, \dots, m-1$, then 
\begin{equation}\label{subcase 1 bi}
b_i= -1 + \sum_{\substack{k >j }}^m a_{ki}, \quad \text{and} \quad b_n=-1.
\end{equation}

 \textbf{Subcase 2:} Let $\det A_{\tilde{\textbf{v}}} \neq (-1)^m$. If $i \neq N_j$ for $j=1, \dots, m$ then $b_i=0$ following Case 2 of the proof of Theorem \ref{Theorem on the matrix A}. Also if $i = N_j$ for $j=1, \dots, m-1$, then 
\begin{equation}\label{subcase 2 bi}
 b_i=\frac{(-1)^m}{\det A_{\tilde{\textbf{v}}}} (a_{(j+1)N_j}-1) \quad \text{and} \quad b_n=\frac{(-1)^m}{\det A_{\tilde{\textbf{v}}}} (a_{1n}-1).
 \end{equation}
\end{proof}

\section{Properties of the cohomology rings of quasitoric manifolds over a vertex cut of a finite product of simplices}\label{Sec_Cohomology of quasitoric manifolds over the vertex cut of a finite product of simplices}

In this section, we calculate the integral cohomology (possibly with minimal generators) of quasitoric manifolds over a vertex cut of a finite product of simplices depending on the determinant of $A_{\tilde{\textbf{v}}}$ defined in \eqref{Eq_matrix A_tilde v}. 
Moreover, we show several relations among the product of the generators. We also classify the integral cohomology rings of quasitoric manifolds over a vertex cut of finite product of simplices.

Let $P=\prod_{j=1}^m \Delta^{n_j}$ be a finite product of simplices as in \eqref{Eq_P is a product of simplices} and $\widebar{P}$ the vertex cut of $P$ at the vertex $\tilde{\textbf{v}}=v_{n_1 n_2 \dots n_m} \in V(P)$. Let $\widebar{\lambda} \colon \mathcal{F}(\widebar{P}) \to \ZZ^n$ be a characteristic function defined as in \eqref{Eq_characteristic function over widebarP} such that $\widebar{\lambda}$ satisfies the following:
\begin{align}
     \det A_{\textbf{{u}}}=\begin{cases}
-1 \quad \text{if } D(\textbf{{u}},\textbf{{u}}_0)=\text{odd}\\
+1 \quad \text{if } D(\textbf{{u}},\textbf{{u}}_0)=\text{even}
\end{cases}
 \end{align}
 where $\textbf{u}, \textbf{u}_0=u_{0 \dots 0} \in V(\widebar{P})$.
 Then by Theorem \ref{Theorem on the matrix A} and Theorem \ref{Theorem on the vector b}, we can classify the matrix $A_{\tilde{\textbf{v}}}$ and the vector $\textbf{b}$ assigned to the new facet $\widebar{F}$. 
Recall the facet set  
$$\mathcal{F}(\widebar{P})=\{\widebar{F}^j_{k_j} := F^j_{k_j} \cap \widebar{P} ~|~ F^j_{k_j} \in \mathcal{F}(P)\} \cup \{\widebar{F}\}$$ of $\widebar{P}$ from \eqref{Eq_vertex and facet set of bar P}.
We assign the indeterminate $x_i$ to the facet $\widebar{F}^j_{k_j}$ where $i=N_{j-1}+k_j$ for $1 \leq k_j \leq n_j$ and $j=1,\dots,m$. 
We also assign the indeterminate $x_{n+i}$ to the facet $F^i_0$.
Moreover, the new facet $\widebar{F}$ of $\widebar{P}$ is assigned to an indeterminate $x$.
Then following \cite{DJ}, we have $$H^*(X(\widebar{P}, \widebar{\lambda})) \cong \ZZ[x_1, \dots,x_n, \dots, x_{n+m},x]/ \widebar{I} + \widebar{J}$$ where $\widebar{I}$ and $\widebar{J}$ are ideals given by the following.
The ideal $\widebar{I}$ can be determined by the minimal non-faces of $\widebar{P}$ 
which are the following 
 \begin{align*}
     &\prod_{k_j=0}^{n_j} \widebar{F}^j_{k_j}   &\text{ for } 1 \leq j \leq m, \\
     &\widebar{F}^j_{n_j} \widebar{F}    &\text{ for } 1 \leq j \leq m, \\
     &\prod_{j=1}^m \prod_{k_j=0}^{n_j-1} \widebar{F}^j_{k_j}.
 \end{align*}
Then the ideal $\widebar{I}$ is generated by the following monomials in $\ZZ[x_1, \dots, x_{n+m},x]$:
\begin{align}
   x_{N_{j-1}+1} \dots x_{N_j} x_{n+j}, \quad
   x x_{N_j}, \quad \text{ and } 
   \prod_{\substack{i=1 \\ i \neq N_j}}^{n+m} x_i
\end{align}
for $j=1, \dots, m$. 
The characteristic matrix is already determined by Proposition \ref{lemma on matrix} and Theorem \ref{Theorem on the matrix A}. Let us consider the vector equation
\begin{equation}\label{Eq_lamda bar j}
\lambda_{\widebar{J}}=\sum_{i=1}^n e_i x_i + \sum_{i=1}^m \textbf{a}_i x_{n+i} + \textbf{b} x.
\end{equation}
Then the ideal $\widebar{J}$ is generated by the coordinates of $\lambda_{\widebar{J}}$. 
Next we compute the integral cohomology $H^*(X(\widebar{P}, \widebar{\lambda}))$ with less number of generators depending on the determinant of the matrix $A_{\tilde{\textbf{v}}}$.

\textbf{Case 1:} Let $\det A_{\tilde{\textbf{v}}} = (-1)^m$. Then $A_{\tilde{\textbf{v}}}$ can be given by \eqref{Eq_matrix when det is 1}, see proof of Theorem \ref{Theorem on the matrix A}.
So, the ideal $\widebar{J}$ is generated by the following degree one polynomials come from \eqref{Eq_lamda bar j}.
\begin{align}
&    x_i-x_{n+j} + \sum_{k>j}^m a_{ki} x_{n+k} +b_ix &\quad & \text{ if } N_{j-1} < i \leq N_j \text{ for } j=1, \dots, m-1 \quad \text{and} \\
&    x_i-x_{n+m} +b_ix &\quad & \text{ if } N_{m-1} < i \leq N_m. \nonumber
\end{align}
This implies
\begin{align}
    x_i=\begin{cases}
    x_{n+j} - \sum_{k>j}^m a_{ki} x_{n+k} - b_ix \quad  &\text{ if } N_{j-1} < i \leq N_j \text{ for } j=1, \dots, m-1, \\
     x_{n+m} - b_ix \quad \hspace{2.2cm} &\text{ if } N_{m-1} < i \leq N_m
     \end{cases}
\end{align}
 in $H^*(X(\widebar{P}$. For simplicity of notation, we denote the generator $x_{n+j} \in H^*(X(\widebar{P}, \widebar{\lambda}))$ by $y_j$ for $j=1, \dots, m$ and $x$ by $y$.
Then 
\begin{equation}\label{Eq_cohomology for case 0}
    H^*(X(\widebar{P}, \widebar{\lambda})) \cong \ZZ[y_1, \dots, y_m,y]/ \widebar{I}
\end{equation}
where the generators of the ideal $\widebar{I}$ have the form of the following polynomials 
\begin{align}\label{Eq_polynomials for quotient when det=1}
    &y_j \prod_{i > N_{j-1}}^{N_j} \bigg(y_j - \sum_{k >j} a_{ki}y_k - b_i y\bigg) \quad & \text{ for } j=1, \dots, m-1, \\
    &y_m \prod_{i>N_{m-1}}^{N_m} (y_m -b_iy), \quad & \nonumber\\
    &y  \bigg(y_j - \sum_{k >j} a_{kN_j}y_k - b_{N_j} y\bigg) \quad & \text{ for } j=1, \dots, m-1, \nonumber\\
   & y (y_m -b_{N_m}y) \quad \text{and} & \nonumber\\
    & \prod_{\substack{i=1\\i \neq N_j}}^{N_{m-1}} \prod_{j=1}^{m-1} \bigg(y_j - \sum_{k >j} a_{ki}y_k - b_i y\bigg) \prod_{i > N_{m-1}}^{N_m-1} (y_m -b_i y). \quad & \nonumber
\end{align}

\textbf{Case 2:} Let $\det A_{\tilde{v}} \neq (-1)^m$. Then $A_{\tilde{\textbf{v}}}$ can be given by \eqref{Eq_matrix when det is not 1}, see the proof of Theorem \ref{Theorem on the matrix A}. Then the ideal $\widebar{J}$ is generated by the following degree one polynomials
\begin{align}\label{Eq_relation between lambda J}
    &x_i -x_{n+j} +b_i x \quad & & \text{ if } N_{j-1} < i < N_j \text{ for } j=1, \dots, m,  \\
    &x_i -x_{n+j} +a_{(j+1)(N_{j})} x_{n+j+1} +b_i x \quad & &\text{ if } i= N_j \text{ for } j=1, \dots, m-1  \quad \text{and} \nonumber \\
    &x_n -x_{n+m}+a_{1n} x_{n+1} +b_nx. \nonumber \quad & 
\end{align}

From the proof of Theorem \ref{Theorem on the vector b}, we have $b_i=0$ if $N_{j-1} < i < N_j$ for $j=1, \dots, m$ when $\det A_{\tilde{\textbf{v}}} \neq 0$; and $b_i$'s are arbitrary when $\det A_{\tilde{\textbf{v}}}=0$. 
For our purpose, we assume $b_i=0$ for $N_{j-1} < i < N_j$ and $j=1, \dots, m$ if $\det A_{\tilde{\textbf{v}}}=0$.
Then from \eqref{Eq_relation between lambda J}, we get the following
\begin{align}
 x_i = \begin{cases}
 x_{n+j} \quad \quad \hspace{3.8cm}  &\text{ if } N_{j-1} < i < N_j, \text{ for } j=1, \dots, m  \\
   x_{n+j} -a_{(j+1)(N_{j})} x_{n+j+1} - b_i x \quad  &\text{ if } i= N_j, \text{ for } j=1, \dots, m-1  \\
   x_{n+m} - a_{1n} x_{n+1} - b_nx.& \end{cases}
\end{align}
For simplicity of notation, here we denote the generator $x_{n+j} \in H^*(X(\widebar{P}, \widebar{\lambda}))$ by $y_j$ for $j=1, \dots, m$ and $x$ by $y$.
Then 
\begin{equation}\label{Eq_cohomology for case 1}
    H^*(X(\widebar{P}, \widebar{\lambda})) \cong \ZZ[y_1, \dots, y_m,y]/ \widebar{I}
\end{equation}
where the generators of the ideal $\widebar{I}$ have the form of the following polynomials 
\begin{align}\label{Eq_polynomials for quotient}
    &y_{j}^{n_j} (y_{j}-a_{j'N_j}y_{j'} -b_{N_j}y) \quad & \text{ for } j=1, \dots, m-1, \\
     &y_{m}^{n_m}( y_{m} - a_{1N_m} y_{1} - b_{N_m}y),\nonumber \\
    &y(y_{j}-a_{j'N_j}y_{j'} -b_{N_j}y) \quad & \text{ for } j=1, \dots, m-1, \nonumber \\
    &y( y_{m} - a_{1N_m} y_{1} - b_ny), \quad \text{and}\nonumber \\
    &\prod_{j=1}^m y_{j}^{n_j}\nonumber
\end{align} with $j'=j+1$.
Note that if $n_j=1$ for $j=1, \dots, m$, then there is no $i$ such that $N_{j-1} <i < N_j=N_{j-1}+1$.

\begin{theorem}\label{Theorem on generators}
The elements $y_1, \dots, y_m, y$ belong to $H^2(X(\widebar{P}, \widebar{\lambda}))$ and satisfy the following.
\begin{enumerate}
    \item $yy_1=yy_2= \dots =yy_m$, and
    \item $y^2=(-1)^{m+1}(\det A_{\tilde{v}})yy_j$ for any j=1, \dots,m.
\end{enumerate}
\end{theorem}

\begin{proof}
The computation of the cohomology ring of a quasitoric manifold following \cite{DJ} gives $y_1, \dots, y_m, y \in H^2(X(\widebar{P}, \widebar{\lambda}))$.
We shall prove $(1)$ and $(2)$ depending on the value of $\det A_{\tilde{\textbf{v}}}$.
There are following two possibilities.

\textbf{Case 1:} Let $\det A_{\tilde{\textbf{v}}}=0$. Following \eqref{Eq_a^j_jk}, \eqref{Eq_matrix when det is not 1} and \eqref{Eq_coordinates are 1}; we get
 \begin{align}\label{Eq_coordinates of c det 0 case}
a_{jk}=
\begin{cases}
-1    \quad & \text{ if } N_{j-1} < k \leq N_j  ~\text{for } j=1, \dots,m , \\
1     \quad & \text{ if } k=N_{j-1} \text{ for } j=2, \dots, n ~\text{or}~ k=n \text{ for } j=1 ,\\
0     \quad & \text{ otherwise}.
\end{cases}
\end{align}
Hence, the polynomials in \eqref{Eq_polynomials for quotient} reduce to the following equations.
\begin{align}\label{Eq_zero case}
y(y_{1}-y_{2}-b_{N_1}y)&=0, \\ \vdots \nonumber\\
y(y_{j}-y_{j'}-b_{N_j}y)&=0, \nonumber \\ \vdots \nonumber \\
y(y_{m-1}-y_{m}-b_{N_{m-1}}y)&=0, \nonumber\\
y(y_{m}-y_{1}-b_{N_m}y)&=0\nonumber
\end{align}
where $j'=j+1$.
By summing up these equations altogether, we get $$(\sum_{j=1}^m b_{N_j}) y^2=0.$$
From \textbf{Case 1} of Theorem \ref{Theorem on the vector b}, we have $\sum_{j=1}^m b_{N_j}=-1$ when $\det A_{\tilde{\textbf{v}}}=0$. This implies $y^2=0$, which proves $(2)$ for this case. Now using $y^2=0$ in the equations in \eqref{Eq_zero case} one by one, we obtain the desired result $(1)$.

\textbf{Case 2:} Let $\det A_{\tilde{\textbf{v}}} \neq 0$. We have following two subcases. 

\textbf{Subcase 1:} Let $\det A_{\tilde{\textbf{v}}} = (-1)^m$. The third and fourth types of polynomials in \eqref{Eq_polynomials for quotient when det=1} become the following $m$ equations
\begin{align}\label{Eq_det 1 case}
    y  \bigg(y_1 - \sum_{k >1} a_{kN_1}y_k - b_{N_1} y\bigg)&=0,\\
    \vdots \nonumber \\
    y  \bigg(y_j - \sum_{k >j} a_{kN_j}y_k - b_{N_j} y\bigg)&=0, \nonumber\\
    \vdots \nonumber \\
    y (y_{m-1}-a_{mN_{m-1}} y_m -b_{N_{m-1}}y)&=0, \nonumber \\
    y (y_m -b_{N_m}y) &=0. \nonumber
\end{align}
We denote the $j$-th equation in \eqref{Eq_det 1 case} by $(\ref{Eq_det 1 case}.j)$ for $j=1, \dots, m$. Recall from the proof of Theorem \ref{Theorem on the vector b} and \eqref{subcase 1 bi} that $b_{N_m}=-1$ for this subcase. Using this in $(\ref{Eq_det 1 case}.m)$ we get $y^2+yy_m=0$. 
Also from \eqref{subcase 1 bi}, we have $b_{N_{m-1}}=-1+a_{mN_{m-1}}.$ Using $b_{N_{m-1}}=-1+a_{mN_{m-1}}$ and $y^2+yy_m=0$ in $(\ref{Eq_det 1 case}.(m-1))$, we get $y^2+yy_{m-1}=0$. 
Repeating similar arguments inductively, we have $$y^2+yy_j=0$$ for $j=1, \dots, m$. This proves our claims $(1)$ and $(2)$ in this subcase.

\textbf{Subcase 2:} Let $\det A_{\tilde{\textbf{v}}} \neq (-1)^m$ and $\det A_{\tilde{\textbf{v}}} \neq 0$. Then third and fourth types of polynomials in \eqref{Eq_polynomials for quotient} become the following $m$ equations
\begin{align}\label{Eq_non zero case}
y(y_{1}-a_{2N_1}y_{2} -b_{N_1}y)&=0,\\ \vdots \nonumber \\
y(y_{j}-a_{j'N_j}y_{j'} -b_{N_j}y)&=0, \nonumber \\ \vdots \nonumber \\
y(y_{m-1}-a_{(m)(N_m)}y_{m} -b_{N_{n-1}}y)&=0, \nonumber\\
y( y_{m} - a_{1N_m} y_{1} - b_{N_m} y)&=0.\nonumber
\end{align}
where $j'=j+1$. We denote the $j$-th equation in \eqref{Eq_non zero case} by (\ref{Eq_non zero case}.$j$) for $j=1, \dots, m$. First, we fix an $j \in \{1, \dots, m\}$ and consider the equation (\ref{Eq_non zero case}.$j''$) for $j''=j-1$.
Then we substitute the value of $yy_{j''}$ from (\ref{Eq_non zero case}.$j''$) into the equation (\ref{Eq_non zero case}.$(j''-1)$).  
So, we have
\begin{equation}\label{alu1}
yy_{j''-1} -a_{j''N_{j''-1}}(a_{jN_{j''}} y_j +b_{N_{j''}} y)y - b_{N_{j''-1}}y^2=0.
\end{equation}
From the proof of Theorem \ref{Theorem on the vector b} and \eqref{subcase 2 bi}, we have $$b_{N_{j''}}=\frac{(-1)^m}{\det A_{\tilde{\textbf{v}}}} (a_{jN_{j''}}-1) \quad \text{ and } \quad b_{N_{j''-1}}=\frac{(-1)^m}{\det A_{\tilde{\textbf{v}}}} (a_{j''N_{j''-1}}-1).$$
Using $b_{N_{j''}}$ and $b_{N_{j''-1}}$ in \eqref{alu1}, we get the following
\begin{align}
& yy_{j-2} - \frac{(a_{j''N_{j''-1}})(a_{jN_{j''}})}{\det A_{\tilde{\textbf{v}}}} [(\det A_{\tilde{\textbf{v}}}) yy_j + (-1)^m y^2] + \frac{(-1)^m}{\det A_{\tilde{\textbf{v}}}}(a_{j''N_{j''-1}})y^2  \\ 
& \hspace{6cm} - \frac{(-1)^m}{\det A_{\tilde{\textbf{v}}}}(a_{j''N_{j''-1}})y^2+ \frac{(-1)^m}{\det A_{\tilde{\textbf{v}}}} y^2=0 \nonumber 
\end{align}
This implies
\begin{equation}\label{alu2}
     yy_{j-2} - \frac{(a_{j''N_{j''-1}})(a_{jN_{j''}})}{\det A_{\tilde{\textbf{v}}}} [(\det A_{\tilde{\textbf{v}}}) yy_j + (-1)^m y^2] + \frac{(-1)^m}{\det A_{\tilde{\textbf{v}}}} y^2=0
\end{equation}
Then we replace the value of $yy_{j-2}$ from \eqref{alu2} in (\ref{Eq_non zero case}.$(j''-2)$). We follow this recursive process of substitution for $(m-1)$ times which gives us the following.
\begin{equation}\label{alu3}
(1-a_{1n}a_{2N_1} \dots a_{mN_{m-1}})[yy_j + \frac{(-1)^m}{\det A_{\tilde{\textbf{v}}}}y^2]=0.
\end{equation}
Note that if $j \neq m$ then in the recursive process at some point we reach (\ref{Eq_non zero case}.$1$), then we plug in the value of $yy_1$ from (\ref{Eq_non zero case}.$1$) in (\ref{Eq_non zero case}.$m$) and carry on the process.
As $\det A_{\tilde{\textbf{v}}}=(-1)^m (1-a_{1n}a_{2N_1} \dots a_{mN_{m-1}}) \neq 0$, from \eqref{alu3} we have 
\begin{align}\label{Eq_next 1}
yy_j + \frac{(-1)^m}{\det A_{\tilde{\textbf{v}}}}y^2=0 \quad \text{ and } \quad
y^2= (-1)^{m+1} (\det A_{\tilde{\textbf{v}}}) yy_j 
\end{align}
for $j=1, \dots, m$; which proves $(2)$ for this subcase.
Also, by \eqref{Eq_next 1} $$yy_1=yy_2= \dots =yy_m.$$
\end{proof}


In the rest of this section, we discuss the relation among the cohomology rings of quasitoric manifolds over a vertex cut of a finite product of simplices. 
Let $X(\widebar{P}, \widebar{\lambda})$ be a quasitoric manifold such that $\det A_{\tilde{\textbf{v}}}=0$ with $b_i=0$ for $i \neq N_j$ and $j=1, \dots, m$ in the vector $\textbf{b}$.
Then by \eqref{Eq_cohomology for case 1} and \eqref{Eq_coordinates of c det 0 case}, the integral cohomology ring of $X(\widebar{P}, \widebar{\lambda})$ is given by
\begin{equation}\label{Eq_bar I next}
H^*(X(\widebar{P}, \widebar{\lambda}); \ZZ) \cong \ZZ[y_1, \dots, y_m,y] / \widebar{I}
\end{equation} 
where the ideal $\widebar{I}$ is generated by 
\begin{align*}
&(1)~ \prod_{j=1}^m y_j^{n_j},\\
&(2)~ y_j^{n_j}(y_j-y_{j+1}-b_{N_j} y),\\
&(3)~ y (y_j-y_{j+1}-b_{N_j} y)
\end{align*}
for $j=1, \dots, m$ with $y_{m+1}:=y_{1}$. Note that the ideal $\widebar{I}$ depends on the entries of the vector $\textbf{b}=\widebar{\lambda}(\widebar{F})$.

\begin{theorem}\label{Th_classify rings when det =0}
Let $P=\prod_{j=1}^m \Delta^{n_j}$ be a product of simplices as in \eqref{Eq_P is a product of simplices} and $\widebar{P}$ is a vertex cut of $P$ along a vertex $\tilde{\emph{\textbf{v}}}=v_{n_1 \dots n_m}$ such that $\det A_{\tilde{\emph{\textbf{v}}}}=0$.
Then the cohomology rings $H^*(X(\widebar{P}, \widebar{\lambda}))$ are isomorphic to each other if $b_i=0$ for $i \neq N_j$ and $j=1, \dots, m$ in the vector $\textbf{b}$ assigned to the new facet $\widebar{F}$.
\end{theorem}

\begin{proof}
Let us define two quasitoric manifolds over $\widebar{P}$ by two characteristic functions $\widebar{\lambda}$ and $\widebar{\lambda}'$.
Let $(A,\textbf{b})$ and $(A', \textbf{b}')$ be the characteristic matrices associated to $\widebar{\lambda}$ and $\widebar{\lambda}'$ respectively.
In Theorem \ref{Theorem on the matrix A}, we have classified the matrices $A$ and $A'$ up to conjugation.
Also from Theorem \ref{Theorem on the vector b}, we get $\sum_{j=1}^m b_{N_j}=-1$; and from the hypothesis, $b_i=0$ for $i \neq N_j$.
Let $\ell $ be an integer which satisfies $\sum_{j=1}^{m-1} (m-j) b_{N_j} \equiv \ell~ (\text{mod}~ m)$ with $0 \leq \ell \leq m-1$.

Let us consider $\textbf{b}'=(b_1', \dots, b_n')^t$ such that $b'_{N_\ell}=-1 \quad \text{and} \quad b_i'= 0 \text{ for } i \neq N_{\ell}.$
Then $\sum_{j=1}^{m-1} (m-j ) b'_{N_j} = \ell$. Thus we have $$\sum_{j=1}^{m-1} (m-j) b_{N_j} \equiv \ell \equiv \sum_{j=1}^{m-1} (m-j ) b'_{N_j} ~ (\text{mod}~ m).$$
This implies that there exists an integer $k$ such that 
\begin{equation}\label{in next eq}
\sum_{j=1}^{m-1} (m-j) b_{N_j} - \sum_{j=1}^{m-1} (m-j ) b'_{N_j} + mk=0.
\end{equation}
Let $\widebar{I}_{\textbf{b}}$ and $\widebar{I}_{\textbf{b}'}$ be the ideals of $\ZZ[y_1, \dots, y_m,y]$ such that $$H^*(X(\widebar{P},\widebar{\lambda}))=\ZZ[y_1, \dots, y_m,y]/\widebar{I}_{\textbf{b}} \quad \text{and} \quad H^*(X(\widebar{P},\widebar{\lambda}'))=\ZZ[y_1, \dots, y_m,y]/\widebar{I}_{\textbf{b}'}$$ where $\widebar{I}_{\textbf{b}}$ and $\widebar{I}_{\textbf{b}'}$ can be described similar to $\widebar{I}$ as in \eqref{Eq_bar I next}.
Now let $$c_j := \sum_{i=1}^j (b_{N_j}-b'_{N_j})+k \text{ for } j=1, \dots, m-1 \quad \text{ and } \quad c_m:=k.$$ So, $c_j$'s are integers for $j=1, \dots,m$ and they satisfy the following two conditions
\begin{align}\label{next eq}
\sum_{j=1}^m c_j=0~(\text{using } \eqref{in next eq}) \quad \text{and}\\
c_j -c_{j-1} = b_{N_j}-b'_{N_j} \text{ for } j=1, \dots, m \nonumber
\end{align}
where $c_0:=c_m$.

Now we consider an automorphism $\phi \colon \ZZ[y_1, \dots, y_m,y] \to \ZZ[y_1, \dots, y_m,y]$ which is defined by $$\phi(y_j)=y_j+c_j y \text{ for } j=1, \dots, m \quad \text{ and } \quad \phi(y)=y.$$
We now show that this $\phi $ induces an isomorphism $$ \widetilde{\phi} \colon H^*(X(\widebar{P},\widebar{\lambda}')) \to H^*(X(\widebar{P},\widebar{\lambda})).$$
For that, we need to show the generators of $\widebar{I}_{\textbf{b}'}$ maps to $\widebar{I}_{\textbf{b}}$ through $\phi$. 
Let us inspect the map $\phi$ on the generators of $\widebar{I}_{\textbf{b}'}$ which can be obtained from \eqref{Eq_bar I next} replacing $b_{N_j}$ by $b'_{N_j}$ one by one. First,
\begin{align*}
\phi(y_j^{n_j}(y_j-y_{j+1}-b'_{N_j} y))&=(y_j+c_j y)^{n_j}(y_j-y_{j+1}+(c_j-c_{j+1}-b'_{N_j})y)\\
&=(y_j+c_j y)^{n_j}(y_j-y_{j+1}-b_{N_j}y) \quad(\text{using } \eqref{next eq})\\
&=0 \quad \text{in } H^*(X(\widebar{P},\widebar{\lambda})).
\end{align*} for $j=1, \dots ,m$.
Similar calculation gives 
\begin{equation*}
\phi(y(y_j-y_{j+1}-b'_{N_j} y))=y(y_j-y_{j+1}-b_{N_j} y)=0 \quad \text{in } H^*(X(\widebar{P},\widebar{\lambda})).
\end{equation*} for $j=1, \dots, m$.
Next we have
\begin{align*}
\phi(\prod_{j=1}^m y_j)&=\prod_{j=1}^m (y_j+c_j y) \\
&=\sum_{j=1}^m (c_jy \prod_{i \neq j}^m y_i) \quad (\text{as }y^2=0, \text{ from Theorem } \ref{Theorem on generators})\\
&=(\sum_{j=1}^m c_j) yy_1^{m-1} \quad (\text{as }yy_j=yy_1 \text{ for any } j, \text{ from Theorem } \ref{Theorem on generators})\\
&=0 \quad (\text{ using } \eqref{next eq}) 
\end{align*}
in $H^*(X(\widebar{P},\widebar{\lambda}))$. Thus we have shown that $\phi$ induces a graded ring homomorphism $\widetilde{\phi} \colon H^*(X(\widebar{P},\widebar{\lambda}')) \to H^*(X(\widebar{P},\widebar{\lambda}))$. Similarly, we can construct the inverse $\widetilde{\phi}^{-1}$ of $\widetilde{\phi}$ such that $$\widetilde{\phi}^{-1}(y_j)=y_j -c_j y \quad \text{and} \quad \widetilde{\phi}^{-1}(y)=y$$ for $j=1, \dots, m$. This proves our claim. 
\end{proof}

\begin{remark}
If $n_j=1$ for $j=1, \dots, m$ in Theorem \ref{Th_classify rings when det =0}, i.e. $P$ is an $m$-cube, one can have a similar result Lemma 5.1 in \cite{HKMP}.
\end{remark}

\begin{theorem}\label{Theorem for additive basis}
Let $P=\prod_{j=1}^m \Delta^{n_j}$ be a finite product of simplices as in \eqref{Eq_P is a product of simplices} and $n_{j_1}, \dots, n_{j_s}$ are greater than $1$ for some $\{j_1, \dots, j_s\} \subseteq \{1, \dots, m\}$. Let $\widebar{P}$ be a vertex cut of $P$ at $\tilde{\textbf{v}}=v_{n_1} \dots v_{n_m}$ and $\widebar{\lambda}$ a characteristic function as considered in Section \ref{Sec_Quasitoric manifolds over vertex cut}. Then the following elements forms an additive basis of $H^4(X(\widebar{P}, \widebar{\lambda}))$:
\begin{enumerate}
    \item $y_iy_j$ where $1 \leq i <j \leq m$, 
    \item $y_j^2$ if $n_j \in \{n_{j_1}, \dots, n_{j_s}\}$, and
    \item $yy_j$ for any $j=1, \dots, m$. 
\end{enumerate}
\end{theorem}

\begin{proof}
From \eqref{Eq_polynomials for quotient} we have that the polynomials $y_iy_j~(1 \leq i <j \leq m)$, $y_j^2$ and $yy_j$ for $j=1, \dots, m$ generate $H^4(X(\widebar{P}, \widebar{\lambda}))$. Theorem \ref{Theorem on generators} implies that $yy_i=yy_j$ for $i,j=1, \dots, m$.
If $n_j=1$, then $y_j^2$ can be represented as a linear combination of $y_iy_j~(1 \leq i <j \leq m)$ and $yy_j$; see \eqref{Eq_polynomials for quotient when det=1} and \eqref{Eq_polynomials for quotient}. So, we do not need to consider $y_j^2$ as generators if $n_j=1$. Thus the number of generators is $\binom{m}{2} +s+1$.

Let us denote the rank of $H^4(X(\widebar{P}, \widebar{\lambda}))$ by $R(H^4(X(\widebar{P}, \widebar{\lambda})))$. We show $$R(H^4(X(\widebar{P}, \widebar{\lambda})))=\binom{m}{2} +s+1$$ to conclude our claim. 
The equivariant connected sum of $\prod_{j=1}^{m}\mathbb{C}P^{n_j}$ and $\mathbb{C}P^n$ at a $T^n$-fixed point with an orientation reversing map is a quasitoric manifold over $\widebar{P}$ and its $i$-th Betti number is $h_i(\widebar{P})$. Here $h_i(\widebar{P})$ is the $i$-th $h$-vector of $\widebar{P}$, see \cite[Section 3]{DJ}. Note that the $2i$-th Betti number of $X(\widebar{P},\widebar{\lambda})$ is $h_i(\widebar{P})$. 
So the fourth Betti number of the equivariant connected sum of $\prod_{j=1}^{m}\mathbb{C}P^{n_j}$ and $\mathbb{C}P^n$ is $\binom{m}{2} +s+1$. Thus the rank of $H^4(X(\widebar{P}, \widebar{\lambda}))$ is $\binom{m}{2} +s+1$.
\end{proof}

For an element $z \in H^2(X(\widebar{P}, \widebar{\lambda}))$, the annihilator of $z$ is defined by $$\text{Ann}(z)=\{w \in H^2(X(\widebar{P}, \widebar{\lambda})) ~|~ zw=0 \text{ in } H^4(X(\widebar{P}, \widebar{\lambda}))\}.$$ 
Since $\{\widebar{F}, \widebar{F}_{n_j}^j\}$ for $j=1, \dots, m$ are non-faces of $\widebar{P}$, then the rank of $\text{Ann}(cy)$ is $m$ for a non-zero constant $c$. The following theorem discusses the converse.

\begin{theorem}\label{Theorem for case 1}
Let $P=\prod_{j=1}^m \Delta^{n_j}$ be a finite product of simplices as in \eqref{Eq_P is a product of simplices} with $m \geq 2$ and $n \geq 3$ and $X(\widebar{P}, \widebar{\lambda})$ is a quasitoric manifold over the vertex cut at $\tilde{\textbf{\emph{v}}}=v_{n_1 \dots n_m}$ of $P$ as discussed previously in Section \ref{Sec_Quasitoric manifolds over vertex cut}. If $\emph{Ann}(z)$ is of rank $m$ for $z \in H^2(X(\widebar{P}, \widebar{\lambda}))$ and $\det A= (-1)^m$, then $z$ is a constant multiple of $y$.
\end{theorem}

\begin{proof}
As $z \in H^2(X(\widebar{P}, \widebar{\lambda}))$, let $z=\sum_{j=1}^m c_j y_j + cy$ for some integers $c_1, \dots, c_m, c$. Let $w \in \text{Ann}(z)$. Then there exists integers $d_1, \dots, d_m,d$ such that $w= \sum_{j=1}^m d_j y_j+dy$ satisfying 
\begin{equation}\label{Eq_Ann}
    (\sum_{j=1}^m c_j y_j + cy)(\sum_{j=1}^m d_j y_j+dy)=0.
\end{equation}

We have $zw \in H^4(X(\widebar{P}, \widebar{\lambda}))$. Also we have the additive basis of $H^4(X(\widebar{P}, \widebar{\lambda}))$ from Theorem \ref{Theorem for additive basis}.
Let $P$ be the product of $s$-many simplices $\Delta^{n_j}$ with $n_j \geq 2$ and the other simplices are closed intervals.
Without loss of generality, we assume the first $(m-s)$-many simplices are closed intervals. Then the basis elements of $H^4(X(\widebar{P}, \widebar{\lambda}))$ are $$yy_1, \quad y_iy_j~(1 \leq i <j \leq m) \quad \text{and} \quad y_j^2 \quad \text{ for } j=(m-s+1), \dots, m.$$ 

\textbf{Case1:}
Let $s=0$. From \eqref{Eq_Ann} and the fact $yy_1=yy_j$ for $j=1, \dots, m$, we have
\begin{equation}
    \sum_{j=1}^m(c_jd_j) y_j^2 + \sum_{\substack{i=1\\i <j}}^{j=m}(c_i d_j + c_j d_i) y_iy_j + (c(\sum_{j=1}^m d_j)+(\sum_{j=1}^m c_j) d-cd) yy_1=0.
\end{equation}
As the linear equation involves only basis elements, all the coefficients are equal to zero. That is 
\begin{align}\label{Eq_Case s=m}
    c_j d_j& =0, \quad & \text{ for } j=1, \dots,m \\
    c_id_j+c_jd_i &=0, \quad & \text{ for } 1\leq i < j \leq m \nonumber\\
    c(\sum_{j=1}^m d_j)+(\sum_{j=1}^m c_j) d-cd&=0. \nonumber
\end{align}
As the rank of $\text{Ann}(z)$ is $m$, there exists coefficient $d_j \neq 0$ for some $j \in \{1, \dots, m\}$. This implies $c_j=0$ from the first set of equations in \eqref{Eq_Case s=m}. Thus from the second set of equations in \eqref{Eq_Case s=m} we conclude $c_j=0$ for all $j=1, \dots,m$. This proves that $z$ is a constant multiple of $y$. 

\textbf{Case 2:} Let $s>0$. Then $P=I^{m-s} \times \prod_{j=1}^s \Delta^{n_j}$ where $I=[0,1]$. 
Thus from \eqref{Eq_polynomials for quotient}, $y_j^2$ can be written as a linear combination of $y_jy_{j+1}$ and $yy_j$ for $j=1, \dots, (m-s)$. We define an index set $$\mathscr{I}\colon = \{(i,j) \colon i=j+1, i=1, \dots, (m-s)\}.$$ Note that the cardinality $|\mathscr{I}|=m-s$.
From \eqref{Eq_Ann} and the fact $yy_1=yy_j$ for $j=1, \dots, m$, we get
\begin{multline}
    \sum_{j=m-s+1}^m (c_jd_j) y_j^2 + \sum_{\substack{1 \leq i <j \leq m \\ (i,j) \not \in \mathscr{I}}}(c_i d_j + c_j d_i) y_iy_j +  \sum_{(i,j) \in \mathscr{I}}(c_i d_j + c_j d_i + a_{jN_i}c_id_i) y_iy_j \\ + (c(\sum_{j=1}^m d_j)+(\sum_{j=1}^m c_j) d-cd + \sum_{j=1}^{m-s} (c_j d_j)b_{N_j}) yy_1=0.
\end{multline}
Then the coefficients of $y_iy_j$ are zeros. We consider the coefficients of $y_iy_j$ for $(i, j) \in \mathscr{I}$, $y_iy_j$ for $j=i+1, ~i=(m-s+1, \dots, m-1)$ and $y_1y_m$. Considering $d_1, d_2, \dots, d_m$ as variables, one can write these linear equations as the following.
\begin{equation*}
\adjustbox{max width=\textwidth}{
    $\begin{pmatrix}
    c_2+c_1a_{2N_1} & c_1 &  0 & \dots & 0 & 0 &0 & \dots & 0 \\
    0 & c_3+c_2a_{3N_2} & c_2 & \dots & 0 & 0 &0 & \dots & 0\\
    \vdots & \vdots & \vdots & \ddots & \ddots & \vdots &\ddots & \vdots\\
    0 & 0& 0 & \dots & c_{s+1}+c_sa_{s+1N_s} & c_s &0 & \dots & 0\\
    0 & 0& 0 & \dots & 0 & c_{s+2} & c_{s+1} & \dots & 0\\
     \vdots & \vdots & \vdots & \ddots & \ddots & \vdots &\ddots & \vdots\\
     c_m & 0& 0 & \dots & 0 & 0 & 0 & \dots & c_1
    \end{pmatrix}$
     $\begin{pmatrix}
    d_1 \\ d_2 \\ \vdots \\ d_s \\ d_{s+1}  \\ \vdots \\d_m
    \end{pmatrix}$
     =$\begin{pmatrix}
    0\\ 0\\ \vdots \\0 \\ 0\\  \vdots \\0
    \end{pmatrix}$.
}
\end{equation*}
Since $\text{Ann}(z)$ is of rank $m$, then the rank of the $m \times m$ coefficient matrix in the above equation is either $0$ or $1$. Thus, we have $c_j=0$ for $j=1, \dots, m$. So, we conclude that $z$ is a non-zero constant multiple of $y$.
\end{proof}

\begin{remark}
If $s=0$ i.e. $P$ involves no interval in its product then we do not need the assumption $\det A_{\tilde{\textbf{v}}}=(-1)^m$ to prove the claim in Theorem \ref{Theorem for case 1}. \qed
\end{remark}

\begin{theorem}\label{Th_classify ring when det=1}
Let $P=\prod_{j=1}^m \Delta^{n_j}$ be a product of simplices as in \eqref{Eq_P is a product of simplices} with $m \geq 2$, $n \geq 3$ and $\widebar{P}$ is a vertex cut of $P$ along a vertex $\tilde{\emph{\textbf{v}}}$ such that $\det A_{\tilde{\emph{\textbf{v}}}}=(-1)^m$.
Then $H^*(X(\widebar{P}, \widebar{\lambda}))$ and $H^*(X(\widebar{P}, \widebar{\lambda}'))$ are isomorphic as graded rings if and only if $H^*(X(P, \lambda))$ and $H^*(X(P, \lambda'))$ are isomorphic as graded rings.
\end{theorem}

\begin{proof}
First, suppose that $H^*(X(P, \lambda))$ and $H^*(X(P, \lambda'))$ be isomorphic as graded rings. Note that $X(\widebar{P}, \widebar{\lambda})$ is a blowup of $X(P, \lambda)$ along $\pi^{-1}(\tilde{\textbf{v}})$ where $\pi$ is the orbit map defined as in \eqref{Eq_define pi}.
This implies $$X(\widebar{P}, \widebar{\lambda})=X(P, \lambda) \sharp \CC P^n$$ where $\sharp$ denotes the equivariant connected sum of two spaces. Also, note that $\CC P^n$ is equivariantly connected with reversed orientation.
Similarly $X(\widebar{P}, \widebar{\lambda}')$ is a blowup of $X(P, \lambda')$ along $\pi^{-1}(\tilde{\textbf{v}})$.
This implies the isomorphism of $H^*(X(\widebar{P}, \widebar{\lambda}))$ and $H^*(X(\widebar{P}, \widebar{\lambda}'))$ as graded rings. Therefore, the `if' part is concluded.

Let $$\widebar{\Phi}: H^*(X(\widebar{P}, \widebar{\lambda})) \to H^*(X(\widebar{P}, \widebar{\lambda}'))$$ be an isomorphism as graded rings.
Recall that $\{y_1, \dots, y_m,y\}$ is an additive basis of $H^2(X(\widebar{P}, \widebar{\lambda}))$ as well as $H^2(X(\widebar{P}, \widebar{\lambda}'))$.
Thus $\widebar{\Phi}$ induces an automorphism $\Phi$ of $\ZZ[y_1, \dots, y_m,y]$ such that $$\Phi(\widebar{I})=\widebar{I}'$$ where $\widebar{I}$ and $\widebar{I}'$ are ideals generated as in \eqref{Eq_polynomials for quotient} for characteristic functions $\widebar{\lambda}$ and $\widebar{\lambda}'$ respectively.
From Theorem \ref{Theorem for case 1}, we have $\Phi(y)=y$ up to sign. Thus $\Phi$ induces an automorphism $\phi$ of $\ZZ[y_1, \dots , y_m]$ satisfying $\phi(I)=I'$. This proves the `only if' part of our claim. 
\end{proof}

\noindent {\bf Acknowledgment.}
 The first author thanks the `ICSR IIT Madras' for a research grants. The second author thanks the IMSc for PDF fellowship.

\bibliographystyle{amsalpha}
\bibliography{ref-Guddu}

\end{document}